# LOGARITHMIC SOBOLEV INEQUALITY FOR ZERO-RANGE DYNAMICS

By Paolo Dai Pra and Gustavo Posta

*Università di Padova and Politecnico di Milano*

We prove that the logarithmic Sobolev constant for zero-range processes in a box of diameter $L$ grows as $L^2$.

**1. Introduction.** Let $\Lambda$ be a cube in $\mathbb{Z}^d$, and $c : \mathbb{N} \to [0, +\infty)$ be a function such that $c(0) = 0$ and $c(n) > 0$ for every $n > 0$. The zero-range process associated to $c(\cdot)$ is a stochastic system of moving particles in $\Lambda$, which evolves according to the following rule: for each site $x \in \Lambda$, containing $\eta_x$ particles, with probability rate $c(\eta_x)$, one particle jumps from $x$ to one of its nearest neighbors chosen with equal probability. Waiting jump times of different sites are independent. This process conserves the total number of particles $N$; for each $N \geq 1$, the zero-range process restricted on configurations with $N$ particles is a finite irreducible Markov chain, whose unique invariant measure $\nu_\Lambda^N$ is proportional to

$$(1.1) \qquad \prod_{x \in \Lambda} \frac{1}{c(\eta_x)!},$$

where

$$c(n)! = \begin{cases} 1, & \text{for } n = 0, \\ c(n)c(n-1)\cdots c(1), & \text{otherwise.} \end{cases}$$

Moreover, the process is reversible with respect to $\nu_\Lambda^N$. Note that in the special case $c(n) = \lambda n$, the process reduces to a system of independent simple symmetric random walks.

Under suitable growth conditions on $c(\cdot)$ (see [7]), the zero range process may be defined in the infinite lattice $\mathbb{Z}^d$. In this case the extremal invariant measures form a one parameter family of infinite product measures, with

marginals

$$\mu_\rho[\eta_x = k] = \frac{1}{Z(\alpha(\rho))} \frac{\alpha(\rho)^k}{c(k)!}, \tag{1.2}$$

where $\rho \geq 0$, $Z(\alpha(\rho))$ is the normalization, and $\alpha(\rho)$ is uniquely determined by the condition $\mu_\rho[\eta_x] = \rho$ (we use here the notation $\mu[f]$ for $\int f \, d\mu$).

Zero-range processes have been extensively studied in terms of their hydrodynamic scaling (see [7]). In proving the hydrodynamic limit for interacting particle systems, estimates on the spectral gap of the generator of the process are often essential. Although such estimates can be avoided for zero-range processes, due to their special structure, the question on how the spectral gap depends on the volume arises naturally. Let $\mathcal{E}_{\nu_\Lambda^N}$ be the Dirichlet form associated to the zero-range process in $\Lambda = [0, L]^d \cap \mathbb{Z}^d$ with $N$ particles. Then the following Poincaré inequality holds for each $f \in L^2(\nu_\Lambda^N)$ (see [8]):

$$\nu_\Lambda^N[f, f] \leq C L^2 \mathcal{E}_{\nu_\Lambda^N}(f, f), \tag{1.3}$$

where $C$ may depend on the dimension $d$ but not on $N$ or $L$. In other terms, (1.3) says that the spectral gap of the Markov generator associated to $\mathcal{E}_{\nu_\Lambda^N}$ shrinks proportionally to $\frac{1}{L^2}$, independently of the number of particles $N$.

Our aim in this paper is to strengthen (1.3) by proving the logarithmic Sobolev inequality: for each $f > 0$,

$$\text{Ent}_{\nu_\Lambda^N}(f) \leq C L^2 \mathcal{E}_{\nu_\Lambda^N}(\sqrt{f}, \sqrt{f}), \tag{1.4}$$

where $\text{Ent}_{\nu_\Lambda^N}(f) = \nu_\Lambda^N[f \log f] - \nu_\Lambda^N[f] \log \nu_\Lambda^N[f]$. Note that, in the case $c(n) = \lambda n$, inequality (1.4) is well known to hold. In that case the system is comprised of independent particles. For a single particle ($N = 1$), inequality (1.4) reduces to the well-known estimates of the logarithmic Sobolev constant of a symmetric random walk. The case of $N > 1$ independent particles follows from the tensor property of the logarithmic Sobolev inequality.

By taking hydrodynamic scaling, it is possible to derive from (1.4) various smoothing properties of the nonlinear semigroup associated with the partial differential equation

$$\partial_t u = \tfrac{1}{2} \Delta(\alpha(u)). \tag{1.5}$$

In recent years, Poincaré and logarithmic Sobolev inequalities have been proved for a number of interacting particle systems, both in presence (see [2, 10]) and in absence (see [12, 15]) of conservation laws, when the possible number of particles allowed for any site is bounded. When the particle number per site is unbounded, very few results are available for conservative dynamics, we can only cite [8] and [14]. In particular, to our knowledge, there is no conservative system with unbounded particle number per site for which the logarithmic Sobolev inequality has been proved to have diffusive scaling. The interest for results in this direction is illustrated in [7], page 423. A common tool in the works concerning dynamics with conservation of the number of particles is the use of ergodic properties of associated *nonconservative dynamics*. Such dynamics are chosen in such a way that the *grand canonical* invariant measure for the conservative dynamics at a given particle density ($\mu_\rho$ in our case) is reversible for the nonconservative one. For example, the nonconservative process corresponding to zero-range is the one for which particle numbers at different sites are independent, and, for a given site $x$, $\eta_x$ evolves according to a birth and death process with death rate $c(\eta_x)$ and constant birth rate $\alpha(\rho)$. Let us denote by $\mathcal{E}_{\Lambda,\rho}^{n.c.}$ the Dirichlet form for this process; the following Poincaré inequality is shown in [8]:

$$(1.6) \qquad \mu_\rho[f,f] \leq C\mathcal{E}_{\Lambda,\rho}^{n.c.}(f,f),$$

where $C$ is independent of both $\rho$ and $L$. In the proof of (1.3), that is by induction on $L$, a result even stronger than (1.6) is used in [8]. Similarly, in the proof of the logarithmic Sobolev inequality for Kawasaki dynamics (see [2]), one relies on the fact that the corresponding nonconservative dynamics (Glauber dynamics) satisfy a logarithmic Sobolev inequality with a constant independent on both density and volume. Thus, following the same principle, the proof of (1.4) would require the logarithmic Sobolev inequality

$$(1.7) \qquad \mathrm{Ent}_{\mu_\rho}(f) \leq C\mathcal{E}_{\Lambda,\rho}^{n.c.}(\sqrt{f},\sqrt{f}).$$

However, under the most commonly used conditions on $c(\cdot)$, this inequality *fails*, in the sense that there is no constant $C$ such that (1.7) holds for every $f > 0$ (even for $\rho$ and $L$ fixed). For instance, we already noticed that, in the case $c(n) = \lambda n$, (1.4) holds, but it is not hard to show that (1.7) does not. Technically speaking, this is one of the key difficulties in the proof of (1.4).

The proof of (1.4) is by induction on $|\Lambda|$. Two schemes are available: the so called *martingale method*, introduced in [10], where, roughly speaking, one point is added to $\Lambda$ at each step of the induction; the *duplication method* used in [2], where the volume of $\Lambda$ is doubled at each step. For Kawasaki dynamics, as well for the proof of (1.3), the two methods are essentially equivalent. This is not the case for (1.4): following the martingale method, a serious difficulty arises at a very early stage of the proof, for a simple reason that is illustrated in Section 3.2. Although the essence of our strategy is borrowed from the one in [2], the lack of (1.7) and unboundedness of particle density have required substantial and nontrivial original work.

The proof turns out to be very long and technical. In a first stage we show that the logarithmic Sobolev inequality holds with a constant independent of $N$, that is,

$$\mathrm{Ent}_{\nu_\Lambda^N}(f) \leq s(L)\mathcal{E}_{\nu_\Lambda^N}(\sqrt{f},\sqrt{f}),$$

with no reasonable control on the growth of $s(L)$. This part, that comes for free for Kawasaki dynamics as well as for any system with bounded density of particles, is contained in [4], and we summarize it here without proofs. In the second stage, that we give here in detail, we set up a second induction in $|\Lambda|$ to get the $L^2$-growth, that is well known to be optimal. This second induction relies on the so-called two blocks-type estimates; the idea of using these estimates in the context of logaritmic Sobolev inequalities is due to Lu and Yau [10]. We have organized the paper in such a way that the strategy of the proof is outlined free of the main technicalities. Thus, after having introduced our notation and stated our main result (Section 2), the whole proof is outlined in Section 3. The proofs of the main estimates used in Section 3 are contained in Sections 4 and 5, except for a relevant technical tool, a comparison inequality between $\nu_\Lambda^N$ and $\mu_{N/|\Lambda|}$, that is proved in Section 6.

**2. Notation and main result.** Throughout this paper, for a given probability space $(\Omega, \mathcal{F}, \mu)$ and $f : \Omega \to \mathbb{R}$ measurable, we use the following notation for mean value and covariance:

$$\mu[f] := \int f \, d\mu, \qquad \mu[f, g] := \mu[(f - \mu[f])(g - \mu[g])]$$

and, for $f \geq 0$,

$$\mathrm{Ent}_\mu(f) := \mu[f \log f] - \mu[f] \log \mu[f],$$

where, by convention, $0 \log 0 = 0$. Similarly, for $\mathcal{G}$ a sub-$\sigma$-field of $\mathcal{F}$, we let $\mu[f|\mathcal{G}]$ to denote the conditional expectation, and

$$\mu[f, g|\mathcal{G}] := \mu[(f - \mu[f|\mathcal{G}])(g - \mu[g|\mathcal{G}])|\mathcal{G}]$$

the conditional covariance.

If $A \subset \Omega$, we denote by $\mathbf{1}(x \in A)$ the indicator function of $A$. If $B \subset A$ is *finite* we will write $B \subset\subset A$. For any $x \in \mathbb{R}$, we will write $\lceil x \rceil := \inf\{n \in \mathbb{Z} : n \geq x\}$.

Let $\Lambda$ be a possibly infinite subset of $\mathbb{Z}^d$, and $\Omega_\Lambda = \mathbb{N}^\Lambda$ be the corresponding configuration space, where $\mathbb{N} = \{0, 1, 2, \dots\}$ is the set of natural numbers. Given a configuration $\eta \in \Omega_\Lambda$ and $x \in \Lambda$, the natural number $\eta_x$ will be referred to as the number of particles at $x$. Moreover, if $\Lambda' \subset \Lambda$, $\eta_{\Lambda'}$ will denote the restriction of $\eta$ to $\Lambda'$. For two elements $\sigma, \xi \in \Omega_\Lambda$, the operations $\sigma \pm \xi$ are defined componentwise (for the difference whenever it returns an element of $\Omega_\Lambda$). In what follows, given $x \in \Lambda$, we make use of the special configuration $\delta^x$, having one particle at $x$ and no other particle. For $f : \Omega_\Lambda \to \mathbb{R}$ and $x, y \in \Lambda$, we let

$$\partial_{xy} f(\eta) := f(\eta - \delta^x + \delta^y) - f(\eta).$$

Consider, at a formal level, the operator

$$\mathcal{L}_\Lambda f(\eta) := \sum_{x \in \Lambda} \sum_{y \sim x} c(\eta_x)\, \partial_{xy} f(\eta), \tag{2.1}$$

where $y \sim x$ means $|x - y| = 1$, and $c \colon \mathbb{N} \to \mathbb{R}_+$ is a function such that $c(0) = 0$ and $\inf\{c(n) \colon n > 0\} > 0$. In the case of $\Lambda$ finite, for each $N \in \mathbb{N} \setminus \{0\}$, $\mathcal{L}_\Lambda$ is the infinitesimal generator of an irreducible Markov chain on the finite state space $\{\eta \in \Omega_\Lambda \colon \bar\eta_\Lambda = N\}$, where

$$\bar\eta_\Lambda := \sum_{x \in \Lambda} \eta_x$$

is the total number of particles in $\Lambda$. The unique stationary measure for this Markov chain is denoted by $\nu_\Lambda^N$ and is given by

$$\nu_\Lambda^N[\{\eta\}] := \frac{1}{Z_\Lambda^N} \prod_{x \in \Lambda} \frac{1}{c(\eta_x)!}, \tag{2.2}$$

where $c(0)! := 1$, $c(k)! := c(1) \cdots c(k)$, for $k > 0$, and $Z_\Lambda^N$ is the normalization factor. The measure $\nu_\Lambda^N$ will be referred to as the *canonical measure*. Note that the system is reversible for $\nu_\Lambda^N$, that is, $\mathcal{L}_\Lambda$ is self-adjoint in $L^2(\nu_\Lambda^N)$ or, equivalently, the *detailed balance condition*

$$c(\eta_x)\nu_\Lambda^N[\{\eta\}] = c(\eta_y + 1)\nu_\Lambda^N[\{\eta - \delta_x + \delta_y\}] \tag{2.3}$$

holds for every $x \in \Lambda$ and $\eta \in \Omega_\Lambda$ such that $\eta_x > 0$.

Consider now a real number $\rho > 0$. By the usual statistical mechanics formalism, we introduce the *grand canonical* probability measure with average density of particles $\rho$:

$$\mu_\rho[\{\eta\}] := \frac{\alpha(\rho)^{\bar\eta_\Lambda}\nu_\Lambda^{\bar\eta_\Lambda}[\{\eta\}]}{\sum_{\xi \in \Omega_\Lambda} \alpha(\rho)^{\bar\xi_\Lambda}\nu_\Lambda^{\bar\xi_\Lambda}[\{\xi\}]} = \frac{1}{Z(\alpha(\rho))} \prod_{x \in \Lambda} \frac{\alpha(\rho)^{\eta_x}}{c(\eta_x)!}, \tag{2.4}$$

where $\alpha(\rho)$ is chosen so that $\mu_\rho[\eta_x] = \rho$, $x \in \Lambda$, and $Z(\alpha(\rho))$ is the corresponding normalization. Clearly, $\mu_\rho$ is a product measure with marginals given by (1.2). Note that $\alpha(\rho)$ is the inverse of the function $\alpha \mapsto \frac{1}{Z(\alpha)} \sum_n \frac{n\alpha^n}{c(n)!}$, that is analytic and strictly increasing in some interval of the form $(0, \alpha^*)$. The inverse function theorem for analytic functions guarantees that $\alpha(\rho)$ is well defined and it is an analytic function of $\rho \in [0, \rho^*)$, for some $\rho^* \in (0, +\infty]$. Under mild additional assumptions (e.g., $\lim_n c(n) = +\infty$, see [7], Chapter 2, Lemma 3.3), $\alpha(\cdot)$ is defined in the whole half line $(0, +\infty)$. The normalization will be sometimes denoted by $Z(\alpha)$ or $Z(\rho)$, depending which parameter needs to be emphasized. Note that $\mathcal{L}_\Lambda$ is self-adjoint in $L^2(\mu_\rho)$ too, since identity (2.3) holds with $\mu_\rho$ in place of $\nu_\Lambda^N$.

When $\Lambda$ is infinite, existence and uniqueness in law of a Feller process generated by $\mathcal{L}_\Lambda$ requires conditions on the rate function $c(\cdot)$. It is shown in [7], Section 2.6, that a sufficient condition is that $c(\cdot)$ is Lipschitz:

CONDITION 2.1 (LG).

$$\sup_{k \in \mathbb{N}} |c(k+1) - c(k)| := a_1 < +\infty.$$

As remarked in [8] for the spectral gap, diffusive scaling of the logarithmic Sobolev inequality requires extra-conditions; in particular, our main result would not hold true in the case $c(k) = c\mathbf{1}(k \in \mathbb{N} \setminus \{0\})$. The following condition, which is the same assumed in [8], is a monotonicity requirement on $c(\cdot)$ that rules out the case just mentioned.

CONDITION 2.2 (M). There exists $k_0 > 0$ and $a_2 > 0$ such that $c(k) - c(j) \geq a_2$ for any $j \in \mathbb{N}$ and $k \geq j + k_0$.

We state here without proof some direct consequences of Conditions 2.1 and 2.2. The proofs of some of them can be found in [8].

PROPOSITION 2.3.

1. *There exists $A_0 > 0$ such that $A_0^{-1}k \leq c(k) \leq A_0 k$ for any $k \in \mathbb{N}$.*
2. *Let $\sigma^2(\rho) := \mu_\rho[\eta_x, \eta_x]$, then*

$$(2.5) \qquad 0 < \inf_{\rho > 0} \frac{\sigma^2(\rho)}{\rho} \leq \sup_{\rho > 0} \frac{\sigma^2(\rho)}{\rho} < +\infty.$$

3. *Let $\alpha(\rho)$ be the function appearing in (2.4); then*

$$(2.6) \qquad 0 < \inf_{\rho > 0} \frac{\alpha(\rho)}{\rho} \leq \sup_{\rho > 0} \frac{\alpha(\rho)}{\rho} < +\infty.$$

For $\Lambda$ infinite, the grand canonical measures $\mu_\rho$ are well defined, and, under Condition 2.1, it can be shown that $\mathcal{L}_\Lambda$ can be extended from cylindrical functions to $L^2(\mu_\rho)$; moreover, the grand canonical measures are all stationary for the system. In the rest of the paper the same symbol $\mu_\rho$ will denote both the grand canonical probability measure on any $\Lambda \subset \mathbb{Z}^d$ and its one-dimensional marginal, which is a probability on $\mathbb{N}$.

In what follows, we choose $\Lambda = [0, L]^d \cap \mathbb{Z}^d$. In order to state our main result, we define the Dirichlet form corresponding to $\mathcal{L}_\Lambda$ and $\nu_\Lambda^N$:

$$(2.7) \quad \mathcal{E}_{\nu_\Lambda^N}(f, g) = -\nu_\Lambda^N[f\mathcal{L}_\Lambda g] = \frac{1}{2} \sum_{x \in \Lambda} \sum_{y \sim x} \nu_\Lambda^N[c(\eta_x)\, \partial_{xy} f(\eta)\, \partial_{xy} g(\eta)].$$

THEOREM 2.1. *Assume that Conditions 2.1 and 2.2 hold. Then there exists a constant $C > 0$ that only depends on $a_1, a_2$ and $d$, such that, for every choice of $N \geq 1$, $L \geq 2$ and $f : \Omega_\Lambda \to \mathbb{R}$, $f \geq 0$, we have*

$$(2.8) \qquad \mathrm{Ent}_{\nu_\Lambda^N}(f) \leq CL^2 \mathcal{E}_{\nu_\Lambda^N}(\sqrt{f}, \sqrt{f}).$$

**3. Outline of the proof.** The proof will be given in one dimension. The extension to higher dimensions follows the analogous extension for the spectral gap, that is given in [7], Appendix 3.3.

**3.1.** *Step* 1: *duplication.* The idea is to prove Theorem 2.1 by induction on $|\Lambda|$. Suppose $|\Lambda| = 2L$, so that $\Lambda = \Lambda_1 \cup \Lambda_2$, $|\Lambda_1| = |\Lambda_2| = L$, where $\Lambda_1, \Lambda_2$ are two disjoint adjacent segments in $\mathbb{Z}$. By a basic identity on the entropy, we have

$$(3.1) \qquad \text{Ent}_{\nu_\Lambda^N}(f) = \nu_\Lambda^N[\text{Ent}_{\nu_\Lambda^N[\cdot|\bar{\eta}_{\Lambda_1}]}(f)] + \text{Ent}_{\nu_\Lambda^N}(\nu_\Lambda^N[f|\bar{\eta}_{\Lambda_1}]).$$

Note that $\nu_\Lambda^N[\cdot|\bar{\eta}_{\Lambda_1}] = \nu_{\Lambda_1}^{\bar{\eta}_{\Lambda_1}} \otimes \nu_{\Lambda_2}^{N-\bar{\eta}_{\Lambda_1}}$. Thus, by the tensor property of the entropy (see [1], Theorem 3.2.2),

$$(3.2) \qquad \nu_\Lambda^N[\text{Ent}_{\nu_\Lambda^N[\cdot|\bar{\eta}_{\Lambda_1}]}(f)] \leq \nu_\Lambda^N\Big[\text{Ent}_{\nu_{\Lambda_1}^{\bar{\eta}_{\Lambda_1}}}(f) + \text{Ent}_{\nu_{\Lambda_2}^{N-\bar{\eta}_{\Lambda_1}}}(f)\Big].$$

Now, let $s(L, N)$ be the maximum of the logarithmic Sobolev constant for the zero-range process in volumes $\Lambda$ with $|\Lambda| \leq L$ and less that $N$ particles, that is, $s(L, N)$ is the smallest constant such that

$$\text{Ent}_{\nu_\Lambda^n}(f) \leq s(L, N)\mathcal{E}_{\nu_\Lambda^n}(\sqrt{f}, \sqrt{f}),$$

for all $f > 0$, $|\Lambda| \leq L$ and $n \leq N$. Then, by (3.2),

$$
\begin{aligned}
(3.3) \qquad &\nu_\Lambda^N[\text{Ent}_{\nu_\Lambda^N[\cdot|\bar{\eta}_{\Lambda_1}]}(f)] \\
&\leq s(L, N)\nu_\Lambda^N\Big[\mathcal{E}_{\nu_{\Lambda_1}^{\bar{\eta}_{\Lambda_1}}}(\sqrt{f}, \sqrt{f}) + \mathcal{E}_{\nu_{\Lambda_2}^{N-\bar{\eta}_{\Lambda_1}}}(\sqrt{f}, \sqrt{f})\Big] \\
&= s(L, N)\mathcal{E}_{\nu_\Lambda^N}(\sqrt{f}, \sqrt{f}).
\end{aligned}
$$

Identity (3.1) and inequality (3.3) suggest estimating $s(L, N)$ by induction on $L$. The hardest thing is to make appropriate estimates on the term $\text{Ent}_{\nu_\Lambda^N}(\nu_\Lambda^N[f|\bar{\eta}_{\Lambda_1}])$. Note that this term is the entropy of a function depending only on the number of particles in $\Lambda_1$.

**3.2.** *Step* 2: *logarithmic Sobolev inequality for the distribution of the number of particles in* $\Lambda_1$. Let

$$\gamma_\Lambda^N(n) = \gamma(n) := \nu_\Lambda^N[\bar{\eta}_{\Lambda_1} = n].$$

$\gamma(\cdot)$ is a probability measure on $\{0, 1, \ldots, N\}$ that is reversible for the birth and death process with generator

$$
\begin{aligned}
(3.4) \qquad \mathcal{A}\varphi(n) := &\Big[\frac{\gamma(n+1)}{\gamma(n)} \wedge 1\Big](\varphi(n+1) - \varphi(n)) \\
&+ \Big[\frac{\gamma(n-1)}{\gamma(n)} \wedge 1\Big](\varphi(n-1) - \varphi(n))
\end{aligned}
$$

and Dirichlet form

$$\mathcal{D}(\varphi, \varphi) = -\langle \varphi, \mathcal{A}\varphi \rangle_{L^2(\gamma)} = \sum_{n=1}^{N} [\gamma(n) \wedge \gamma(n-1)](\varphi(n) - \varphi(n-1))^2.$$

Logarithmic Sobolev inequalities for birth and death processes are studied in [13]. The nontrivial proof that conditions in [13] are satisfied by $\gamma(n)$ leads to the following result, whose proof is in [4].

PROPOSITION 3.1. *The Markov chain with generator* (3.4) *has a logarithmic Sobolev constant proportional to* $N$, *that is, there exists a constant* $C > 0$ *such that, for all* $\varphi \geq 0$,

$$\text{Ent}_\gamma(\varphi) \leq CN\mathcal{D}(\sqrt{\varphi}, \sqrt{\varphi}).$$

REMARK 3.2. At this point it is possible to see why we use the duplication approach of [2] instead of the martingale method of Lu and Yau [10]. Indeed, assume $\Lambda_2$ is as above, $\Lambda_1$ is a single point and $c(n) = \lambda n$. In this case $\gamma(\cdot)$ is a Binomial with $N$ trials and probability of success $1/L$. It is well known (see [5]) that the associated birth and death process (3.4) has a logarithmic Sobolev constant of order $N \log(L)$. The extra factor $\log(L)$ does not get fixed by the rest of our argument.

We now apply Proposition 3.1 to the second summand of the right-hand side of (3.1), and we obtain

$$\text{Ent}_{\nu_\Lambda^N}(\nu_\Lambda^N[f|\bar{\eta}_{\Lambda_1}])$$

$$(3.5) \quad \begin{aligned} &\leq CN \sum_{n=1}^{N} [\gamma(n) \wedge \gamma(n-1)] \Big[ \sqrt{\nu_\Lambda^N[f|\bar{\eta}_{\Lambda_1} = n]} - \sqrt{\nu_\Lambda^N[f|\bar{\eta}_{\Lambda_1} = n-1]} \Big]^2 \\ &\leq CN \sum_{n=1}^{N} \frac{\gamma(n) \wedge \gamma(n-1)}{\nu_\Lambda^N[f|\bar{\eta}_{\Lambda_1} = n] \vee \nu_\Lambda^N[f|\bar{\eta}_{\Lambda_1} = n-1]} \\ &\qquad \times (\nu_\Lambda^N[f|\bar{\eta}_{\Lambda_1} = n] - \nu_\Lambda^N[f|\bar{\eta}_{\Lambda_1} = n-1])^2, \end{aligned}$$

where we have used the inequality $(\sqrt{x} - \sqrt{y})^2 \leq \frac{(x-y)^2}{x \vee y}$, $x, y > 0$.

3.3. *Step* 3: *study of the term* $\nu_\Lambda^N[f|\bar{\eta}_{\Lambda_1} = n] - \nu_\Lambda^N[f|\bar{\eta}_{\Lambda_1} = n-1]$. One of the key points in the proof of Theorem 2.1 consists in finding the "right" representation for the discrete gradient $\nu_\Lambda^N[f|\bar{\eta}_{\Lambda_1} = n] - \nu_\Lambda^N[f|\bar{\eta}_{\Lambda_1} = n-1]$, which appears in the right-hand side of (3.5). The following result is proved in [4].

PROPOSITION 3.3. *For every $f$ and every $n = 1, 2, \ldots, N$, we have*

$$\nu_\Lambda^N[f|\bar\eta_{\Lambda_1} = n] - \nu_\Lambda^N[f|\bar\eta_{\Lambda_1} = n - 1]$$

$$\text{(3.6)} \qquad = \frac{\gamma(n-1)}{\gamma(n)} \frac{1}{nL} \left( \nu_\Lambda^N\left[ \sum_{x \in \Lambda_1, y \in \Lambda_2} h(\eta_x) c(\eta_y) \, \partial_{yx} f \Big| \bar\eta_{\Lambda_1} = n-1 \right] \right.$$

$$\left. + \nu_\Lambda^N\left[ f, \sum_{x \in \Lambda_1, y \in \Lambda_2} h(\eta_x) c(\eta_y) \Big| \bar\eta_{\Lambda_1} = n-1 \right] \right),$$

*where*

$$h(n) := \frac{n+1}{c(n+1)}.$$

*Moreover, by exchanging the roles of $\Lambda_1$ and $\Lambda_2$, the right-hand side of* (3.6) *can be equivalently written as, for every $n = 0, 1, \ldots, N-1$,*

$$-\frac{\gamma(N-n)}{\gamma(N-n+1)} \frac{1}{(N-n+1)L}$$

$$\text{(3.7)} \qquad \times \left( \nu_\Lambda^N\left[ \sum_{x \in \Lambda_1, y \in \Lambda_2} h(\eta_y) c(\eta_x) \, \partial_{xy} f \Big| \bar\eta_{\Lambda_1} = n \right] \right.$$

$$\left. + \nu_\Lambda^N\left[ f, \sum_{x \in \Lambda_1, y \in \Lambda_2} h(\eta_y) c(\eta_x) \Big| \bar\eta_{\Lambda_1} = n \right] \right).$$

The representations (3.6) and (3.7) will be used for $n \geq \frac{N}{2}$ and $n < \frac{N}{2}$, respectively. For convenience, we rewrite (3.6) and (3.7) as

$$\text{(3.8)} \qquad \nu_\Lambda^N[f|\bar\eta_{\Lambda_1} = n] - \nu_\Lambda^N[f|\bar\eta_{\Lambda_1} = n-1] =: A(n) + B(n),$$

where

$$\text{(3.9)} \; A(n) := \begin{cases} \dfrac{\gamma(n-1)}{\gamma(n)} \dfrac{1}{nL} \\ \times \nu_\Lambda^N\left[ \displaystyle\sum_{x \in \Lambda_1, y \in \Lambda_2} h(\eta_x) c(\eta_y) \, \partial_{yx} f \Big| \bar\eta_{\Lambda_1} = n-1 \right], & \text{for } n \geq \dfrac{N}{2}, \\[1em] -\dfrac{\gamma(N-n)}{\gamma(N-n+1)} \dfrac{1}{(N-n+1)L} \\ \times \nu_\Lambda^N\left[ \displaystyle\sum_{x \in \Lambda_1, y \in \Lambda_2} h(\eta_y) c(\eta_x) \, \partial_{xy} f \Big| \bar\eta_{\Lambda_1} = n \right], & \text{for } n < \dfrac{N}{2}, \end{cases}$$

and

$$(3.10) \ B(n) := \begin{cases} \dfrac{\gamma(n-1)}{\gamma(n)} \dfrac{1}{nL} \\ \quad \times \nu_\Lambda^N\Big[f, \displaystyle\sum_{x \in \Lambda_1, y \in \Lambda_2} h(\eta_x)c(\eta_y) \Big| \bar\eta_{\Lambda_1} = n-1\Big], \quad \text{for } n \geq \dfrac{N}{2}, \\ -\dfrac{\gamma(N-n)}{\gamma(N-n+1)} \dfrac{1}{(N-n+1)L} \\ \quad \times \nu_\Lambda^N\Big[f, \displaystyle\sum_{x \in \Lambda_1, y \in \Lambda_2} h(\eta_y)c(\eta_x) \Big| \bar\eta_{\Lambda_1} = n\Big], \quad \text{for } n < \dfrac{N}{2}. \end{cases}$$

Thus, our next aim is to get estimates on the two terms in the right-hand side of (3.8). It is useful to stress that the two terms are qualitatively different. The first term, $A(n)$, contains discrete gradients of $f$. It is mainly this term that is responsible for the growth $L^2$ of the logarithmic Sobolev constant. Estimates on $A(n)$ are essentially insensitive to the precise form of $c(\cdot)$. Indeed, the dependence of $A(n)$ on $L$ and $N$ is of the same order as in the case $c(n) = \lambda n$, that is, the case of independent particles. Quite differently, the term $B(n)$ vanishes in the case of independent particles, since, in that case, the term $\sum_{x \in \Lambda_1, y \in \Lambda_2} h(\eta_x)c(\eta_y)$ is a.s. constant with respect to $\nu_\Lambda^N[\cdot|\bar\eta_{\Lambda_1} = n-1]$. Thus, $B(n)$ somewhat depends on interaction between particles. Note that our model is not necessarily a "small perturbation" of a system of independent particles; there is no small parameter in the model that guarantees that $B(n)$ is small enough. Essentially all technical results of this paper are concerned with estimating $B(n)$.

3.4. *Step* 4: *estimates on $A(n)$.*  The following proposition, proved in [4], gives the key estimate on $A(n)$.

PROPOSITION 3.4.  *There is a constant $C > 0$ such that*

$$A^2(n) \leq \frac{CL^2}{N}(\nu_\Lambda^N[f|\bar\eta_{\Lambda_1} = n] \vee \nu_\Lambda^N[f|\bar\eta_{\Lambda_1} = n-1])$$

$$\times \left[\frac{\gamma(n-1)}{\gamma(n)}\mathcal{E}_{\nu_\Lambda^N[\cdot|\bar\eta_{\Lambda_1} = n-1]}(\sqrt{f}, \sqrt{f}) + \mathcal{E}_{\nu_\Lambda^N[\cdot|\bar\eta_{\Lambda_1} = n]}(\sqrt{f}, \sqrt{f})\right].$$

Let us try to see where we are now. Let us ignore, for the moment, the term $B(n)$, that is, let us pretend that $B(n) \equiv 0$. Thus, by (3.8) and Proposition 3.4, we would have

$$(\nu_\Lambda^N[f|\bar\eta_{\Lambda_1} = n] - \nu_\Lambda^N[f|\bar\eta_{\Lambda_1} = n-1])^2$$

$$(3.11) \qquad \leq \frac{CL^2}{N}(\nu_\Lambda^N[f|\bar\eta_{\Lambda_1} = n] \vee \nu_\Lambda^N[f|\bar\eta_{\Lambda_1} = n-1])$$

$$\times \left[ \frac{\gamma(n-1)}{\gamma(n)} \mathcal{E}_{\nu_\Lambda^N[\cdot|\bar\eta_{\Lambda_1}=n-1]}(\sqrt f, \sqrt f) + \mathcal{E}_{\nu_\Lambda^N[\cdot|\bar\eta_{\Lambda_1}=n]}(\sqrt f, \sqrt f) \right].$$

Inserting (3.11) into (3.5), we get, for some possibly different constant $C$,

$$(3.12) \qquad \operatorname{Ent}_{\nu_\Lambda^N}(\nu_\Lambda^N(f|\bar\eta_{\Lambda_1})) \le CL^2 \mathcal{E}_{\nu_\Lambda^N}(\sqrt f, \sqrt f),$$

where we have used the obvious identity:

$$(3.13) \qquad \nu_\Lambda^N[\mathcal{E}_{\nu_\Lambda^N[\cdot|\mathcal{G}]}(\sqrt f, \sqrt f)] = \mathcal{E}_{\nu_\Lambda^N}(\sqrt f, \sqrt f),$$

which holds for any $\sigma$-field $\mathcal{G}$. Inequality (3.12), together with (3.1) and (3.3), yields

$$(3.14) \qquad s(2L, N) \le s(L, N) + CL^2.$$

Thus, if we can show that

$$(3.15) \qquad \sup_N s(2, N) < +\infty,$$

(see Proposition 3.7 next) then Theorem 2.1 would follow from (3.14). In all this, however, we have totally ignored the contribution of $B(n)$.

### 3.5. *Step* 5: *preliminary analysis of* $B(n)$.

We confine ourselves to the analysis of $B(n)$ for $n \ge \frac{N}{2}$, since the case $n < \frac{N}{2}$ is identical. Consider the covariance that appears in the right-hand side of the first formula of (3.10). By elementary properties of the covariance and the fact that $\nu_\Lambda^N[\cdot|\bar\eta_{\Lambda_1}] = \nu_{\Lambda_1}^{\bar\eta_{\Lambda_1}} \otimes \nu_{\Lambda_2}^{N-\bar\eta_{\Lambda_1}}$, we get

$$
\nu_\Lambda^N\left[ f, \sum_{x\in\Lambda_1, y\in\Lambda_2} h(\eta_x)c(\eta_y)\Big|\bar\eta_{\Lambda_1}=n-1 \right]
$$

$$(3.16) \qquad = \nu_{\Lambda_1}^{n-1}\left[ \sum_{x\in\Lambda_1} h(\eta_x)\nu_{\Lambda_2}^{N-n+1}\left[ f, \sum_{y\in\Lambda_2} c(\eta_y) \right] \right]$$

$$+ \nu_{\Lambda_1}^{n-1}\left[ \nu_{\Lambda_2}^{N-n+1}[f], \sum_{x\in\Lambda_1} h(\eta_x) \right] \nu_{\Lambda_2}^{N-n+1}\left[ \sum_{y\in\Lambda_2} c(\eta_y) \right].$$

It follows by Conditions 2.1 and 2.2 (see Proposition 2.3) that, for some constant $C > 0$, $h(\eta_x) \le C$ and $c(\eta_y) \le C\eta_y$. Thus, a simple estimate on the two summands in (3.16) yields, for some $C > 0$,

$$
B^2(n) \le \frac{\gamma^2(n-1)}{\gamma^2(n)} \left( \frac{C}{n^2} \nu_{\Lambda_2}^{N-n+1}\left[ \nu_{\Lambda_1}^{n-1}[f], \sum_{y\in\Lambda_2} c(\eta_y) \right]^2 \right.
$$

$$(3.17) \qquad \left. + \frac{C(N-n+1)^2}{n^2 L^2} \nu_{\Lambda_1}^{n-1}\left[ \nu_{\Lambda_2}^{N-n+1}[f], \sum_{x\in\Lambda_1} h(\eta_x) \right]^2 \right).$$

Thus, our next aim is to estimate the two covariances in (3.17). We establish two levels of estimates for these covariances. We first recall the "rough estimates" on $B(n)$ obtained in [4], which use elementary properties of $\nu_\Lambda^N$. These estimates, together with the duplication argument mentioned above, will allow to prove that, for all $L$,

$$(3.18) \qquad\qquad \sup_N s(L, N) \equiv s(L) < +\infty,$$

without, however, any reasonable information on the growth of $s(L)$ in $L$. Then we will obtain sharper estimates (two-blocks estimates) on $B(n)$ in the case of $L$ large enough, say, $L \geq L_0$. Repeating the duplication argument starting from size $L_0$, we will then succeed in obtaining the right $L^2$-growth of $s(L)$. Note that, to begin this duplication argument, we need to know that $\sup_N s(N, L_0) < +\infty$, which is by no means a trivial fact. It should be remarked that, for dynamics with exclusion rules [2, 16], this problem does not appear.

3.6. *Rough estimates on $B(n)$: entropy inequality and estimates on moment generating functions.* By (3.17), estimating $B^2(n)$ consists in estimating two covariances. In general, covariances can be estimated by the following *entropy inequality* that holds for every probability measure $\nu$ (see [1], Section 1.2.2):

$$(3.19) \qquad \nu[f, g] = \nu[f(g - \nu[g])] \leq \frac{\nu[f]}{t} \log \nu[e^{t(g-\nu[g])}] + \frac{1}{t} \operatorname{Ent}_\nu(f),$$

where $f \geq 0$ and $t > 0$ is arbitrary. Since in (3.17) we need to estimate the square of a covariance, we write (3.19) with $-g$ in place of $g$, and obtain

$$(3.20) \quad |\nu[f, g]| \leq \frac{\nu[f]}{t} \log(\nu[e^{t(g-\nu[g])}] \vee \nu[e^{-t(g-\nu[g])}]) + \frac{1}{t} \operatorname{Ent}_\nu(f).$$

Therefore, we first get estimates on the moment generating functions $\nu[e^{\pm t(g-\nu[g])}]$, and then optimize (3.20) over $t > 0$.

In [2], for Kawasaki dynamics, estimates for moment generating functions with respect to $\nu_\Lambda^N$ are obtained in two main steps: 1. One replaces the canonical ensemble $\nu_\Lambda^N$ with the corresponding grand canonical ensemble, by using results on equivalence of ensembles. 2. The nonconservative (Heat Bath) dynamics for the grand canonical ensemble have a logarithmic Sobolev constant that is uniform in the volume, from which estimates on moment generating functions follow (see [1], Section 7.4.1) by the so-called Herbst argument. In our model there are difficulties for both steps. For step 1, we can prove equivalence of ensemble only for $\Lambda$ large, so the case of small $\Lambda$ must be treated separately. For step 2, as we pointed out in the Introduction, there is no logarithmic Sobolev inequality for the nonconservative dynamics corresponding to the zero-range process. Actually, it would not be hard to show that good estimates for moment generating functions follow, via Herbst argument, from a weaker *modified logarithmic Sobolev inequality* (the entropy inequality in [3]). This inequality holds true [3] in the case $c(n) = \lambda n$, but nothing is known under our assumptions on $c(\cdot)$.

Our approach is to give estimates on *special* moment generating functions, where ad-hoc arguments can be used. Note that the covariances in (3.17) involve functions of $\eta_{\Lambda_1}$ or $\eta_{\Lambda_2}$. In the next two propositions we write $\Lambda$ for $\Lambda_1$ and $\Lambda_2$, and denote by $N$ the number of particles in $\Lambda$. Their proof can be found in [4].

PROPOSITION 3.5. *Let $x \in \Lambda$. Then there is a constant $A_L$ depending on $L = |\Lambda|$ such that, for every $N > 0$ and $t \in [-1, 1]$,*

$$(3.21) \qquad \nu_\Lambda^N[e^{t(c(\eta_x) - \nu_\Lambda^N[c(\eta_x)])}] \le e^{A_L N t^2}$$

*and*

$$(3.22) \qquad \nu_\Lambda^N[e^{tN(h(\eta_x) - \nu_\Lambda^N[h(\eta_x)])}] \le A_L e^{A_L(Nt^2 + \sqrt{N}|t|)}.$$

Using (3.19), (3.20) and (3.22) and optimizing over $t > 0$, we get estimates on the covariances appearing in (3.17).

PROPOSITION 3.6. *There exists a constant $C_L$ depending on $L = |\Lambda|$ such that the following inequalities hold:*

$$(3.23) \qquad \nu_\Lambda^N\left[f, \sum_{x \in \Lambda} c(\eta_x)\right]^2 \le C_L N \nu_\Lambda^N[f] \operatorname{Ent}_{\nu_\Lambda^N}(f),$$

$$(3.24) \qquad \nu_\Lambda^N\left[f, \sum_{x \in \Lambda} h(\eta_x)\right]^2 \le C_L N^{-1} \nu_\Lambda^N[f][\nu_\Lambda^N[f] + \operatorname{Ent}_{\nu_\Lambda^N}(f)].$$

Inserting these new estimates in (3.17), we obtain, for some possibly different $C_L$,

$$
\begin{aligned}
B^2(n) \le{}& \frac{C_L \gamma^2(n-1)}{\gamma^2(n)} \nu_\Lambda^N[f|\bar\eta_{\Lambda_1} = n-1] \\
(3.25) \qquad & \times \Bigg( \frac{N-n+1}{n^2} \nu_{\Lambda_1}^{n-1}[\operatorname{Ent}_{\nu_{\Lambda_2}^{N-n+1}}(f)] \\
& + \frac{(N-n+1)^2}{n^3} (\nu_\Lambda^N[f|\bar\eta_{\Lambda_1} = n-1] + \nu_{\Lambda_2}^{N-n+1}[\operatorname{Ent}_{\nu_{\Lambda_1}^n}(f)]) \Bigg).
\end{aligned}
$$

In order to simplify (3.25), we use Proposition 7.5 in [4], which gives

$$(3.26) \qquad \frac{n}{C(N-n+1)} \le \frac{\gamma(n-1)}{\gamma(n)} \le \frac{Cn}{N-n+1}$$

for some $C > 0$. It follows that

$$\frac{\gamma^2(n-1)}{\gamma^2(n)}\frac{N-n+1}{n^2} \leq \frac{C^2}{N-n+1} \leq \frac{\gamma(n-1)}{\gamma(n)}\frac{C^3}{n}$$

and

$$\frac{\gamma^2(n-1)}{\gamma^2(n)}\frac{(N-n+1)^2}{n^3} \leq \frac{C^2}{n}.$$

Thus, (3.25) implies, for some $C_L > 0$ depending on $L$, recalling also that $n \geq \frac{N}{2}$,

$$N\frac{\gamma(n) \wedge \gamma(n-1)}{\nu_\Lambda^N[f|\bar\eta_{\Lambda_1} = n] \vee \nu_\Lambda^N[f|\bar\eta_{\Lambda_1} = n-1]}B^2(n)$$

$$(3.27) \quad \leq C_L\gamma(n-1)(\nu_\Lambda^N[f|\bar\eta_{\Lambda_1} = n-1]$$

$$+ \nu_{\Lambda_1}^{n-1}[\mathrm{Ent}_{\nu_{\Lambda_2}^{N-n+1}}(f)] + \nu_{\Lambda_2}^{N-n+1}[\mathrm{Ent}_{\nu_{\Lambda_1}^n}(f)]).$$

Now, we bound the two terms

$$\mathrm{Ent}_{\nu_{\Lambda_1}^{n-1}}(f) \quad \text{and} \quad \mathrm{Ent}_{\nu_{\Lambda_2}^{N-n+1}}(f)$$

by, respectively,

$$s(N,L)\mathcal{E}_{\nu_{\Lambda_1}^{n-1}}(\sqrt f, \sqrt f) \quad \text{and} \quad s(N,L)\mathcal{E}_{\nu_{\Lambda_2}^{N-n+1}}(\sqrt f, \sqrt f),$$

and insert these estimates in (3.6). What comes out is then used to estimate (3.5), after having obtained the corresponding estimates for $n < \frac{N}{2}$. Recalling the estimates for $A(n)$, straightforward computations yield

$$\mathrm{Ent}_{\nu_\Lambda^N}(\nu_\Lambda^N[f|\bar\eta_{\Lambda_1}])$$

$$(3.28) \quad \leq C_L\mathcal{E}_{\nu_\Lambda^N}(\sqrt f, \sqrt f) + C_L\nu_\Lambda^N[f] + C_Ls(N,L)\mathcal{E}_{\nu_\Lambda^N}(\sqrt f, \sqrt f),$$

for some $C_L > 0$ depending on $L$. Inequality (3.28), together with (3.3), gives, for a possibly different constant $C_L > 0$,

$$\mathrm{Ent}_{\nu_\Lambda^N}(f)$$

$$(3.29) \quad \leq C_L\mathcal{E}_{\nu_\Lambda^N}(\sqrt f, \sqrt f) + C_L\nu_\Lambda^N[f] + C_Ls(N,L)\mathcal{E}_{\nu_\Lambda^N}(\sqrt f, \sqrt f).$$

To deal with the term $\nu_\Lambda^N[f]$ in (3.29), we use the following well-known argument. Set $\bar f = (\sqrt f - \nu_\Lambda^N[\sqrt f])^2$. By the Rothaus inequality (see [1], Lemma 4.3.8),

$$\mathrm{Ent}_{\nu_\Lambda^N}(f) \leq \mathrm{Ent}_{\nu_\Lambda^N}(\bar f) + 2\nu_\Lambda^N[\sqrt f, \sqrt f].$$

Using this inequality and replacing $f$ by $\bar{f}$ in (3.29), we get, for a different $C_L$,

$$\operatorname{Ent}_{\nu_\Lambda^N}(f)$$

$$
\begin{aligned}
(3.30) \quad &\leq C_L \mathcal{E}_{\nu_\Lambda^N}(\sqrt{f}, \sqrt{f}) + C_L \nu_\Lambda^N[\sqrt{f}, \sqrt{f}] + C_L s(N, L) \mathcal{E}_{\nu_\Lambda^N}(\sqrt{f}, \sqrt{f}) \\
&\leq D_L \mathcal{E}_{\nu_\Lambda^N}(\sqrt{f}, \sqrt{f}) + D_L s(N, L) \mathcal{E}_{\nu_\Lambda^N}(\sqrt{f}, \sqrt{f}),
\end{aligned}
$$

where, in the last line, we have used the Poincaré inequality (1.3). Therefore,

$$(3.31) \qquad s(2L, N) \leq D_L[s(L, N) + 1],$$

that implies (3.18), provided we prove the following "basis step" for the induction.

PROPOSITION 3.7.

$$\sup_N s(2, N) < +\infty.$$

The proof of Proposition 3.7 is also given in [4], Proposition 3.6. As we pointed out above, (3.31) gives no indication on how $s(L) = \sup_N s(N, L)$ grows with $L$.

3.7. *Two blocks estimates on $B(n)$.* Our goal here is to improve the estimates in Proposition 3.6, for large $\Lambda$. Following the terminology in [10, 16], we refer to these estimates as the "two blocks estimates."

We begin by considering the covariance in (3.23). First note that this covariance is left unchanged if we replace $c(n)$ by $\tilde{c}(n) := c(n) - \delta n$, where $\delta$ is a parameter that may depend on $N$ and $L$, and whose value, for reasons that will be apparent later, will be set to $\delta = \alpha'(N/L)$ for $N/L$ not too small, and $\delta = \frac{\alpha(N/L)}{N/L}$ for $N/L$ small. Then we partition $\Lambda$ into disjoint blocks $C_1, C_2, \ldots, C_m$ that we assume, with no loss of generality, to have equal size $l = L/m$. Denote by $\mathcal{G}$ the $\sigma$-field generated by $\bar{\eta}_{C_1}, \bar{\eta}_{C_2}, \ldots, \bar{\eta}_{C_m}$. Note that

$$(3.32) \qquad \nu_\Lambda^N[\cdot | \mathcal{G}] = \nu_{C_1}^{\bar{\eta}_{C_1}} \otimes \nu_{C_2}^{\bar{\eta}_{C_2}} \otimes \cdots \otimes \nu_{C_m}^{\bar{\eta}_{C_m}}.$$

Thus, the usual formula for conditional covariances gives

$$\nu_\Lambda^N[f, g] = \nu_\Lambda^N[\nu_\Lambda^N[f, g | \mathcal{G}]] + \nu_\Lambda^N[f, \nu_\Lambda^N[g | \mathcal{G}]],$$

and we obtain

$$(3.33) \qquad \nu_\Lambda^N\left[f, \sum_{x \in \Lambda} \tilde{c}(\eta_x)\right] = \nu_\Lambda^N\left[\nu_\Lambda^N\left[f, \sum_{x \in \Lambda} \tilde{c}(\eta_x)\Big| \mathcal{G}\right]\right]$$

$$(3.34) \qquad\qquad + \nu_\Lambda^N\left[f, \sum_{k=1}^m \nu_{C_k}^{\bar{\eta}_{C_k}}\left[\sum_{x \in C_k} \tilde{c}(\eta_x)\right]\right].$$

We now treat the two summands (3.33) and (3.34). The term (3.33) is much simpler, and it is actually insensitive to the choice of $\delta$ and $l$. The following result will be proved in Section 4.

PROPOSITION 3.8.   *There is a constant $C$, possibly depending on $l$, such that*

$$(3.35) \qquad \nu_\Lambda^N \left[ \nu_\Lambda^N \left[ f, \sum_{x \in \Lambda} \tilde{c}(\eta_x) \Big| \mathcal{G} \right] \right]^2 \leq C N \nu_\Lambda^N[f] \nu_\Lambda^N [\mathrm{Ent}_{\nu_\Lambda^N[\cdot|\mathcal{G}]}(f)].$$

Now, by (3.18), the zero-range process on $C_k$ with $n_k$ particles satisfies the logarithmic Sobolev inequality with a constant $s$ depending on $l$ but *not* on $n_k$. Moreover, by (3.32) and the tensor property of the entropy (3.2), it follows that

$$\mathrm{Ent}_{\nu_\Lambda^N[\cdot|\mathcal{G}]}(f) \leq \sum_{k=1}^m \nu_\Lambda^N \left[ \mathrm{Ent}_{\nu_{C_k}^{\bar{n}_{C_k}}}(f) | \mathcal{G} \right] \leq s \sum_{k=1}^m \nu_\Lambda^N \left[ \mathcal{E}_{\nu_{C_k}^{\bar{n}_{C_k}}}(\sqrt{f}, \sqrt{f}) | \mathcal{G} \right].$$

Thus, averaging with respect to $\nu_\Lambda^N$, we easily obtain

$$\nu_\Lambda^N [\mathrm{Ent}_{\nu_\Lambda^N[\cdot|\mathcal{G}]}(f)] \leq s \mathcal{E}_{\nu_\Lambda^N}(\sqrt{f}, \sqrt{f}),$$

and from Proposition 3.8, we obtain the following corollary.

COROLLARY 3.9.   *There is a constant $C$, possibly depending on $l$, such that*

$$(3.36) \qquad \nu_\Lambda^N \left[ \nu_\Lambda^N \left[ f, \sum_{x \in \Lambda} \tilde{c}(\eta_x) \Big| \mathcal{G} \right] \right]^2 \leq C N \nu_\Lambda^N[f] \mathcal{E}_{\nu_\Lambda^N}(\sqrt{f}, \sqrt{f}).$$

Estimating (3.34) is much harder. Here the choice $\delta$ becomes essential, as well as the possibility of choosing $l$ and $L$ large. The starting point is again to estimate the covariance in (3.34) by the entropy inequality (3.20). The challenge now is to get sharp estimates on the moment generating function. Both these estimates and the optimization of (3.20) over $t$ require considerable work, that we do no summarize in this preliminary section. The outcome is the following inequality.

PROPOSITION 3.10.   *For every $\varepsilon > 0$, we can choose $l = l_\varepsilon$, $L_0 = L_{0,\varepsilon}$, and a finite constant $C_\varepsilon$ such that, for $L \geq L_0$,*

$$\nu_\Lambda^N \left[ f, \sum_{k=1}^m \nu_{C_k}^{\bar{n}_{C_k}} \left[ \sum_{x \in C_k} \tilde{c}(\eta_x) \right] \right]^2$$
$$\leq N \nu_\Lambda^N[f] [C_\varepsilon \nu_\Lambda^N[f] + C_\varepsilon L^2 \mathcal{E}_{\nu_\Lambda^N}(\sqrt{f}, \sqrt{f}) + \varepsilon \, \mathrm{Ent}_{\nu_\Lambda^N}(f)].$$

By putting together the estimates in Corollary 3.9 and Proposition 3.10, we get the following result.

COROLLARY 3.11.   *For every $\varepsilon > 0$, there exists a finite constant $C_\varepsilon$ such that*

$$
\begin{aligned}
(3.37) \quad & \nu_\Lambda^N \left[ f, \sum_{x \in \Lambda} c(\eta_x) \right]^2 \\
& \leq N \nu_\Lambda^N[f][C_\varepsilon \nu_\Lambda^N[f] + C_\varepsilon L^2 \mathcal{E}_{\nu_\Lambda^N}(\sqrt{f}, \sqrt{f}) + \varepsilon \operatorname{Ent}_{\nu_\Lambda^N}(f)].
\end{aligned}
$$

To complete our two blocks estimate on $B(n)$, we have to estimate the second covariance in (3.17). The argument follows the same lines as for the first covariance.

PROPOSITION 3.12.   *For every $\varepsilon > 0$, there exists a finite constant $C_\varepsilon$ such that*

$$
\begin{aligned}
(3.38) \quad & \nu_\Lambda^N \left[ f, \sum_{x \in \Lambda} h(\eta_x) \right]^2 \\
& \leq \frac{L^2}{N} \nu_\Lambda^N[f][C_\varepsilon \nu_\Lambda^N[f] + C_\varepsilon L^2 \mathcal{E}_{\nu_\Lambda^N}(\sqrt{f}, \sqrt{f}) + \varepsilon \operatorname{Ent}_{\nu_\Lambda^N}(f)].
\end{aligned}
$$

3.8. *Concluding the proof.* Now, as we did in (3.25)–(3.27), we use the inequalities obtained in Corollary 3.11 and Proposition 3.12 in order to estimate $B^2(n)$ (for $n \geq \frac{N}{2}$, the other case being identical), and we obtain the following refinement of (3.27), holding for $L$ sufficiently large:

$$
\begin{aligned}
(3.39) \quad & N \frac{\gamma(n) \wedge \gamma(n-1)}{\nu_\Lambda^N[f | \bar{\eta}_{\Lambda_1} = n] \vee \nu_\Lambda^N[f | \bar{\eta}_{\Lambda_1} = n-1]} B^2(n) \\
& \leq \gamma(n-1)(C_\varepsilon \nu_\Lambda^N[f | \bar{\eta}_{\Lambda_1} = n-1] \\
& \qquad\qquad + C_\varepsilon L^2 \mathcal{E}_{\nu_\Lambda^N[\cdot | \bar{\eta}_{\Lambda_1} = n-1]}(\sqrt{f}, \sqrt{f}) \\
& \qquad\qquad + \varepsilon \operatorname{Ent}_{\nu_\Lambda^N[\cdot | \bar{\eta}_{\Lambda_1} = n-1]}(f)).
\end{aligned}
$$

Thus, inserting these estimates in (3.6) and then in (3.5), we obtain, for some $C > 0$ independent of $\varepsilon$,

$$
\begin{aligned}
(3.40) \quad & \operatorname{Ent}_{\nu_\Lambda^N}(\nu_\Lambda^N[f | \bar{\eta}_{\Lambda_1}]) \\
& \leq CL^2 \mathcal{E}_{\nu_\Lambda^N}(\sqrt{f}, \sqrt{f}) + C_\varepsilon \nu_\Lambda^N[f] \\
& \qquad + C_\varepsilon L^2 \mathcal{E}_{\nu_\Lambda^N}(\sqrt{f}, \sqrt{f}) + C\varepsilon \operatorname{Ent}_{\nu_\Lambda^N}(\nu_\Lambda^N[f | \bar{\eta}_{\Lambda_1}]).
\end{aligned}
$$

Thus, by (3.1), (3.3) and (3.40), choosing $\varepsilon$ so that $C\varepsilon < 1$, eliminating the term $\nu_\Lambda^N[f]$ as in (3.30), for some $C > 0$, we get

$$(3.41) \qquad \mathrm{Ent}_{\nu_\Lambda^N}(f) \leq s(L,N)\mathcal{E}_{\nu_\Lambda^N}(\sqrt{f}, \sqrt{f}) + CL^2 \mathcal{E}_{\nu_\Lambda^N}(\sqrt{f}, \sqrt{f}),$$

that gives

$$s(2L, N) \leq s(L, N) + CL^2,$$

from which the conclusion follows.

The rest of this paper is devoted to the proof of all results concerning two-blocks estimates that have been previously stated without proof.

## 4. Two blocks estimates I: proofs of Propositions 3.8 and 3.10.

We begin with the proof of Proposition 3.8, which is a simple modification of the one of Proposition 3.6, given in [4].

PROOF OF PROPOSITION 3.8. By (3.32) and the entropy inequality (3.19), we have

$$\nu_\Lambda^N \left[ f, \sum_{x \in \Lambda} c(\eta_x) \Big| \mathcal{G} \right]$$

$$\leq \frac{\nu_\Lambda^N[f|\mathcal{G}]}{t} \sum_{k=1}^m \log \nu_{C_k}^{\bar{\eta}_{C_k}} \left[ \exp\left( t \sum_{x \in C_k} \left( c(\eta_x) - \nu_{C_k}^{\bar{\eta}_{C_k}}[c(\eta_x)] \right) \right) \right]$$

$$+ \frac{1}{t} \mathrm{Ent}_{\nu_\Lambda^N[\cdot|\mathcal{G}]}(f).$$

We now estimate each term $\nu_{C_k}^{\bar{\eta}_{C_k}}[\exp(t\sum_{x \in C_k}(c(\eta_x) - \nu_{C_k}^{\bar{\eta}_{C_k}}[c(\eta_x)]))]$ using Proposition 3.5, obtaining, for $0 < t \leq 1$,

$$\nu_\Lambda^N \left[ f, \sum_{x \in \Lambda} c(\eta_x) \Big| \mathcal{G} \right] \leq \frac{\nu_\Lambda^N[f|\mathcal{G}]}{t} CNt^2 + \frac{1}{t}\mathrm{Ent}_{\nu_\Lambda^N[\cdot|\mathcal{G}]}(f),$$

for some constant $C$ depending on $l$ but not on $N$. After averaging over $\nu_\Lambda^N$, we get

$$(4.1) \qquad \nu_\Lambda^N \left[ \nu_\Lambda^N \left[ f, \sum_{x \in \Lambda} c(\eta_x) \Big| \mathcal{G} \right] \right] \leq CNt\nu_\Lambda^N[f] + \nu_\Lambda^N \left[ \frac{1}{t}\mathrm{Ent}_{\nu_\Lambda^N[\cdot|\mathcal{G}]}(f) \right].$$

At this point we observe that (4.1) holds with $-c$ in place of $c$. Thus, in this last inequality we can square both sides, and then optimize over $t$ as follows. Set

$$t_*^2 = \frac{\nu_\Lambda^N[\mathrm{Ent}_{\nu_\Lambda^N[\cdot|\mathcal{G}]}(f)]}{N\nu_\Lambda^N[f]}.$$

If $t_* \leq 1$, then plugging $t_*$ in (4.1), we get the inequality in (3.35), for a possibly different $C > 0$. For $t_* > 1$, using point 1 of Proposition 2.3, we get the simple bounds, for some $C > 0$,

$$\nu_\Lambda^N \left[ \nu_\Lambda^N \left[ f, \sum_{x \in \Lambda} c(\eta_x) \Big| \mathcal{G} \right] \right]^2 \leq C N^2 \nu_\Lambda^N [f]^2 \leq C N^2 t_*^2 \nu_\Lambda^N [f]^2$$

$$= C N \nu_\Lambda^N [f] \nu_\Lambda^N [\mathrm{Ent}_{\nu_\Lambda^N [\cdot | \mathcal{G}]}(f)],$$

and (3.35) follows in this case too. $\square$

We now turn to the proof of Proposition 3.10. In this proof we set $\rho = \frac{N}{L}$. As customary in two block estimates [2, 8], we treat separately the case of small density.

The following lemma will be used in various parts of the proof.

LEMMA 4.1. *For every $\rho > 0$ and $t \in \mathbb{R}$,*

$$\mu_\rho[e^{t(c-\alpha)}] \leq e^{\alpha a_1 t^2 e^{a_1 |t|}},$$

*where $\alpha = \alpha(\rho)$ and $a_1$ is the constant appearing in Condition 2.1.*

PROOF. We follow the well-known Herbst argument. Let $T$ be the shift defined by $Tf(n) = f(n+1)$. A simple computation shows that

$$(4.2) \qquad\qquad \mu_\rho[cf] = \alpha \mu_\rho[Tf].$$

Moreover, by Condition 2.1,

$$(4.3) \qquad\qquad \|Tc - c\|_{+\infty} \leq a_1.$$

Now let

$$\phi(t) := \mu_\rho[e^{tc}].$$

Suppose, first, $t > 0$. We have

$$t\phi'(t) - \phi(t)\log\phi(t) = t\mu_\rho[ce^{tc}] - \mu_\rho[e^{tc}]\log\mu_\rho[e^{tc}]$$

$$\leq \alpha t \mu_\rho[e^{t(Tc)} - e^{tc}],$$

where we have used (4.2) and the fact that, by Jensen's inequality,

$$-\log\mu_\rho[e^{tc}] \leq -\alpha t.$$

Thus, using the inequality $|e^x - e^y| \leq |x - y|e^{|x-y|}e^x$ and (4.5), we get

$$t\phi'(t) - \phi(t)\log\phi(t) \leq \alpha a_1 t^2 e^{a_1 t}\phi(t).$$

Letting $\psi(t) := \frac{\log \phi(t)}{t}$, this last inequality becomes

$$\psi'(t) \le \alpha a_1 e^{a_1 t}.$$

Since

$$\lim_{t \downarrow 0} \psi(t) = \mu_\rho(c) = \alpha,$$

we have

$$\psi(t) \le \alpha e^{a_1 t}$$

that becomes

$$\mu_\rho[e^{t(c-\alpha)}] = \phi(t) e^{-\alpha t} \le e^{\alpha t(e^{a_1 t} - 1)} \le e^{\alpha a_1 t^2 e^{a_1 t}}.$$

To deal with a negative $t$, we need to show that, for $t > 0$,

$$\mu_\rho[e^{-t(c-\alpha)}] \le e^{\alpha a_1 t^2 e^{a_1 t}},$$

which is shown by repeating the above argument. $\quad\square$

REMARK 4.2. The estimates on $\phi(t)$ in the proof of Lemma 4.1 and the fact that $c(n) \ge \varepsilon n$ for some $\varepsilon > 0$ imply that, for some $C > 0$,

$$\mu_\rho[e^{\eta_x}] \le e^{Ct\rho e^{Ct}}.$$

4.1. *Case of density away from zero.* In this section we assume $\rho \ge \rho_0 > 0$, where $\rho_0$ is a constant that will be fixed later. Some of the constants in the estimates that follow may depend on $\rho_0$. In what follows we assume $\nu_\Lambda^N[f] = 1$. We begin with the entropy inequality, for $t > 0$:

$$\nu_\Lambda^N\left[f, \sum_{k=1}^m \nu_{C_k}^{\bar{\eta}_{C_k}}\left[\sum_{x \in C_k} \tilde{c}(\eta_x)\right]\right]$$

$$(4.4) \quad \le \frac{1}{t} \log \nu_\Lambda^N\left[\exp\left(t \sum_{k=1}^m \left(\nu_{C_k}^{\bar{\eta}_{C_k}}\left[\sum_{x \in C_k} \tilde{c}(\eta_x)\right] - \nu_\Lambda^N\left[\sum_{x \in C_k} \tilde{c}(\eta_x)\right]\right)\right)\right]$$

$$+ \frac{1}{t} \text{Ent}_{\nu_\Lambda^N}(f).$$

We have to estimate a moment generating function. One of the ingredients to improve the estimate given in Proposition 3.5 is to replace the canonical ensemble $\nu_\Lambda^N$ by the grand-canonical ensemble $\mu_\rho$. We argue as follows. Define

$$a_k := \nu_{C_k}^{\bar{\eta}_{C_k}}\left[\sum_{x \in C_k} \tilde{c}(\eta_x)\right] - \nu_\Lambda^N\left[\sum_{x \in C_k} \tilde{c}(\eta_x)\right].$$

By the Cauchy–Schwarz inequality

$$\nu_\Lambda^N\left[\exp\left(t\sum_{k=1}^m\left(\nu_{C_k}^{\bar\eta_{C_k}}\left[\sum_{x\in C_k}\tilde c(\eta_x)\right]-\nu_\Lambda^N\left[\sum_{x\in C_k}\tilde c(\eta_x)\right]\right)\right)\right]$$

$$(4.5)\qquad =\nu_\Lambda^N\left[\exp\left(t\sum_{k=1}^m a_k\right)\right]$$

$$\le\left\{\nu_\Lambda^N\left[\exp\left(2t\sum_{1\le k<m/2}a_k\right)\right]\nu_\Lambda^N\left[\exp\left(2t\sum_{m/2\le k<m}a_k\right)\right]\right\}^{1/2}.$$

Now, for both factors in (4.5), assuming $L$ large enough, we may use the inequality between ensembles given in the following proposition, whose proof is left for the final section of this paper. For generality, this proposition is stated (and proved) in any dimension $d\ge 1$.

PROPOSITION 4.1. *Fix $\delta_0\in(0,1)$; then there exist $v_0$ and $A_0>0$ such that, for any $\Lambda'\subset\Lambda\subset\subset\mathbb{Z}^d$ with $|\Lambda|\ge v_0$ and $|\Lambda'|/|\Lambda|\le\delta_0$,*

$$(4.6)\qquad \nu_\Lambda^N[\eta_{\Lambda'}=\sigma_{\Lambda'}]\le A_0\mu_{N/|\Lambda|}[\eta_{\Lambda'}=\sigma_{\Lambda'}],$$

*for any $\sigma_{\Lambda'}\in\Omega_{\Lambda'}$, $N\in\mathbb{N}\setminus\{0\}$.*

Thus, using (4.6) with $\rho=N/L$ in (4.5), we obtain, for some $C>0$,

$$\nu_\Lambda^N\left[\exp\left(t\sum_{k=1}^m\left(\nu_{C_k}^{\bar\eta_{C_k}}\left[\sum_{x\in C_k}\tilde c(\eta_x)\right]-\nu_\Lambda^N\left[\sum_{x\in C_k}\tilde c(\eta_x)\right]\right)\right)\right]$$

$$(4.7)\qquad \le C\left\{\mu_\rho\left[\exp\left(2t\sum_{1\le k<m/2}a_k\right)\right]\mu_\rho\left[\exp\left(2t\sum_{m/2\le k<m}a_k\right)\right]\right\}^{1/2}.$$

The next step consists in replacing, in the terms $a_k$ in (4.7), the expectation $\nu_\Lambda^N[\sum_{x\in C_k}\tilde c(\eta_x)]$ by the corresponding $\mu_\rho$-expectation. We recall that, according to Corollary 6.4 in [8], there is a constant $C>0$ such that

$$(4.8)\qquad |\nu_\Lambda^N[c(\eta_x)]-\mu_\rho[c(\eta_x)]|\le\frac{C}{L}\sqrt{1+\rho}.$$

Thus, using (4.8) in (4.7) and the fact that $\mu_\rho$ is a product measure, we get

$$\nu_\Lambda^N\left[\exp\left(t\sum_{k=1}^m\left(\nu_{C_k}^{\bar\eta_{C_k}}\left[\sum_{x\in C_k}\tilde c(\eta_x)\right]-\nu_\Lambda^N\left[\sum_{x\in C_k}\tilde c(\eta_x)\right]\right)\right)\right]$$

$$(4.9)\qquad \le Ce^{Ct\sqrt N}\mu_\rho\left[\exp\left(2t\left(\nu_{C_1}^{\bar\eta_{C_1}}\left[\sum_{x\in C_1}\tilde c(\eta_x)\right]-\mu_\rho\left[\sum_{x\in C_1}\tilde c(\eta_x)\right]\right)\right)\right]^{L/l}.$$

To get estimates on the expectation in (4.9), observe that, adding and subtracting $\mu_{\bar{\eta}_{C_1}/l}[\sum_{x \in C_1} \tilde{c}(\eta_x)]$ and using the Cauchy–Schwarz inequality, we obtain

$$
\begin{aligned}
(4.10) \quad & \mu_\rho\left[\exp\left\{t\left(\nu_{C_1}^{\bar{\eta}_{C_1}}\left[\sum_{x \in C_1} \tilde{c}(\eta_x)\right] - \mu_\rho\left[\sum_{x \in C_1} \tilde{c}(\eta_x)\right]\right)\right\}\right] \\
& \leq \{\mu_\rho[e^{2tY}]\mu_\rho[e^{2tW}]\}^{1/2},
\end{aligned}
$$

where

$$
\begin{aligned}
Y := {}& \nu_{C_1}^{\bar{\eta}_{C_1}}\left[\sum_{x \in C_1} \tilde{c}(\eta_x)\right] - \mu_{\bar{\eta}_{C_1}/l}\left[\sum_{x \in C_1} \tilde{c}(\eta_x)\right] \\
& - \mu_\rho\left[\nu_{C_1}^{\bar{\eta}_{C_1}}\left[\sum_{x \in C_1} \tilde{c}(\eta_x)\right] - \mu_{\bar{\eta}_{C_1}/l}\left[\sum_{x \in C_1} \tilde{c}(\eta_x)\right]\right]
\end{aligned}
$$

and

$$
\begin{aligned}
W := {}& \mu_{\bar{\eta}_{C_1}/l}\left[\sum_{x \in C_1} \tilde{c}(\eta_x)\right] - \mu_\rho\left[\sum_{x \in C_1} \tilde{c}(\eta_x)\right] \\
& - \mu_\rho\left[\mu_{\bar{\eta}_{C_1}/l}\left[\sum_{x \in C_1} \tilde{c}(\eta_x)\right] - \mu_\rho\left[\sum_{x \in C_1} \tilde{c}(\eta_x)\right]\right].
\end{aligned}
$$

Now our aim is to get separate estimates on the two factors $\mu_\rho[e^{2tY}]$ and $\mu_\rho[e^{2tW}]$ in (4.10).

4.1.1. *Estimates of the factor* $\mu_\rho[e^{2tY}]$ *in* (4.10). In this section we will show that, for a suitable constant $C > 0$, the first factor of (4.10) is bounded above by

$$
(4.11) \qquad \exp[C\rho t^2 e^{Ct}\sqrt{\rho}e^{Ct}].
$$

Let $X$ be a real random variable, and $E(X)$ its expectation. By a straightforward Taylor expansion, the following inequality holds:

$$
(4.12) \qquad E(e^X) \leq \exp[E(X) + \tfrac{1}{2}E(X^2 e^{|X|})].
$$

Using (4.12) and (4.8), we obtain

$$
\begin{aligned}
(4.13) \quad & \mu_\rho[e^{tY}] \leq \exp\left\{Ct^2\mu_\rho\left[\left(\nu_{C_1}^{\bar{\eta}_{C_1}}\left[\sum_{x \in C_1} \tilde{c}(\eta_x)\right] - \mu_{\bar{\eta}_{C_1}/l}\left[\sum_{x \in C_1} \tilde{c}(\eta_x)\right]\right)^2 e^{tY}\right]\right\} \\
& \leq \exp\left\{Ct^2\mu_\rho\left[\left(1 + \frac{\bar{\eta}_{C_1}}{l}\right)e^{Ct\sqrt{1+\bar{\eta}_{C_1}/l}}\right]\right\}
\end{aligned}
$$

$$\leq \exp\left\{ Ct^2 \mu_\rho \left[\left(1 + \frac{\bar{\eta}_{C_1}}{l}\right)^2\right]^{1/2} \mu_\rho [e^{2Ct\sqrt{1+\bar{\eta}_{C_1}/l}}]^{1/2}\right\}$$

$$\leq \exp\{Ct^2 \rho \mu_\rho [e^{2Ct\sqrt{1+\bar{\eta}_{C_1}/l}}]\},$$

where $C > 0$ is some constant that may change in the different steps (we have used here the simple inequality $\mu_\rho[(1 + \frac{\bar{\eta}_{C_1}}{l})^2]^{1/2} \leq C\rho$ for $\rho \geq \rho_0$, which follows from the variance estimates in Proposition 2.3).

The next step consists in estimating the expectation $\mu_\rho[e^{2Ct\sqrt{1+\bar{\eta}_{C_1}/l}}]$ in (4.13). We make use of the following simple lemma.

LEMMA 4.3. *Let $X \geq 0$ be a random variable, and $g: [0, +\infty) \to [0, +\infty)$ be a function such that, for all $t \geq 0$,*

$$(4.14) \qquad E(e^{tX}) \leq e^{tg(t)E(X)}.$$

*Then, for all $t \geq 0$,*

$$(4.15) \qquad E(e^{t\sqrt{X}}) \leq \exp[t\sqrt{2g(2t) + g(t)}\sqrt{E(X)}] + e^t.$$

PROOF. Let $k = 2g(2t) + g(t)$. Then, using the inequality $\sqrt{x} \leq x + 1$, the Cauchy–Schwarz inequality and the Chebyshev inequality,

$$E(e^{t\sqrt{X}}) \leq E(e^{t\sqrt{X}}\mathbf{1}(X < kE(X))) + E(e^{t\sqrt{X}}\mathbf{1}(X \geq kE(X)))$$

$$\leq e^{t\sqrt{kE(X)}} + e^t E(e^{tX}\mathbf{1}(X \geq kE(X)))$$

$$\leq e^{t\sqrt{kE(X)}} + e^t [E(e^{2tX})]^{1/2}[P(X \geq kE(X))]^{1/2}$$

$$\leq e^{t\sqrt{kE(X)}} + e^t e^{tg(2t)E(X)}\sqrt{\frac{E(e^{tX})}{e^{tkE(X)}}}$$

$$\leq e^{t\sqrt{kE(X)}} + e^t \exp\left[\frac{t}{2}E(X)(2g(2t) + g(t) - k)\right]$$

$$\leq e^{t\sqrt{kE(X)}} + e^t. \qquad \square$$

We now go back to the estimate of (4.13). First, by the inequality $\sqrt{1+x} \leq 1 + \sqrt{x}$, we have $\mu_\rho[\exp(2Ct\sqrt{1 + \bar{\eta}_{C_1}/l})] \leq \exp(2Ct)\mu_\rho[\exp(2Ct\sqrt{\bar{\eta}_{C_1}/l})]$. Now, we apply Lemma 4.3 to $X = \bar{\eta}_{C_1}/l$. Note that $\mu_\rho[\bar{\eta}_{C_1}/l] = \rho$ and, by what was observed in Remark 4.2, $\mu_\rho[\exp(t\rho\bar{\eta}_{C_1}/l)] \leq \exp(D\rho t e^{Dt})$, for some $D > 0$. It follows by Lemma 4.3 that

$$\mu_\rho[e^{2Ct\sqrt{\bar{\eta}_{C_1}/l}}] \leq \exp[2C\sqrt{2De^{4CDt} + De^{2CDt}}\sqrt{\rho}t] + e^{2Ct}.$$

Putting it all together, we have that, for a suitable constant $C > 0$, the first factor of (4.10) is bounded above by (4.11).

4.1.2. *Estimates of the factor $\mu_\rho[e^{2tW}]$ in* (4.10). In this section we will prove that, for a suitable constant,

$$(4.16) \qquad \mu_\rho[e^{tW}] \le \exp[Ct^2\rho e^{Ctl}\sqrt{\rho}\, e^{Ct^2l^2\rho e^{Ctl}}].$$

We begin by recalling that $\tilde{c}(n) = c(n) - \delta n$. It is only at this point that the choice of $\delta$ becomes relevant. Note that, since $\mu_\rho[\tilde{c}] = \alpha(\rho) - \delta\rho$, if we define

$$F := \alpha\left(\frac{\bar{\eta}_{C_1}}{l}\right) - \alpha(\rho) - \delta\left(\frac{\bar{\eta}_{C_1}}{l} - \rho\right) - \mu_\rho\left[\alpha\left(\frac{\bar{\eta}_{C_1}}{l}\right) - \alpha(\rho) - \delta\left(\frac{\bar{\eta}_{C_1}}{l} - \rho\right)\right],$$

we have

$$(4.17) \qquad \mu_\rho[e^{tW}] = \mu_\rho[e^{ltF}].$$

Now we set $\delta = \alpha'(\rho)$, in order for $F$ to be a first-order Taylor expansion. To get estimates on $\mu_\rho[e^{ltF}]$, we proceed as follows. Using the inequality $e^x \le 1 + x + \frac{1}{2}x^2 e^{|x|}$ on $e^{ltF}$, and the fact that $\mu_\rho[F] = 0$, we have that

$$\mu_\rho[e^{ltF}] \le 1 + \frac{1}{2}l^2t^2\mu_\rho[F^2 e^{tl|F|}].$$

Then using the inequality $1 + x \le e^x$, which holds for $x > 0$, and the Cauchy–Schwarz inequality, we get

$$(4.18) \qquad \log\mu_\rho[e^{ltF}] \le \frac{1}{2}t^2l^2(\mu_\rho[F^4])^{1/2}(\mu_\rho[e^{2tl|F|}])^{1/2}.$$

We now estimate separately the two expectations in the right-hand side of (4.18). For the estimate of $\mu_\rho[F^4]$, we proceed as in Lemma 3.4 in [8]. We remark that the same expectation was estimated in [8], but their estimate is not good for us. We would rather follow here the deeper argument they used to estimate $\mu_\rho[F^2]$. Let

$$(4.19) \qquad \mathcal{Z} := \frac{1}{l}\sum_{x \in C_1}\frac{\eta_x - \rho}{\sigma(\rho)},$$

where $\sigma^2(\rho) = \mu_\rho[n, n]$. By Lemma 5.2 in [8], each $\mu_\rho$-moment of $\frac{\eta_x - \rho}{\sigma(\rho)}$ is bounded in $\rho \ge \rho_0$. Fix a constant $\beta > 0$.

- For $|\mathcal{Z}| > \beta$, the following rough estimate on $F$ that comes from Lipschitz continuity of $\alpha$ suffices:

$$|F| \le C\sigma(\rho)|\mathcal{Z}|,$$

for some $C > 0$. Thus, possibly modifying the value of $C$, we have

$$\mu_\rho[F^4\mathbf{1}(|\mathcal{Z}| > \beta)]$$
$$\le C\sigma^4(\rho)\mu_\rho[\mathcal{Z}^4\mathbf{1}(|\mathcal{Z}| > \beta)]$$

$$\leq \frac{\sigma^4(\rho)}{\beta^4} \mu_\rho[\mathcal{Z}^8]$$

$$(4.20)$$

$$\leq \frac{C\sigma^4(\rho)}{l^8\beta^4} \left\{ \sum_{x,y,z,w \in C_1} \mu_\rho\left[ \frac{(\eta_x - \rho)^2}{\sigma^2(\rho)} \frac{(\eta_y - \rho)^2}{\sigma^2(\rho)} \frac{(\eta_z - \rho)^2}{\sigma^2(\rho)} \frac{(\eta_w - \rho)^2}{\sigma^2(\rho)} \right] \right.$$

$$\left. + \sum_{x,y,z \in C_1} \mu_\rho\left[ \frac{(\eta_x - \rho)^3}{\sigma^3(\rho)} \frac{(\eta_y - \rho)^3}{\sigma^3(\rho)} \frac{(\eta_z - \rho)^2}{\sigma^2(\rho)} \right] \right\}$$

$$\leq \frac{C(\beta)\rho^2}{l^4},$$

where $C(\beta) > 0$ depends on $\beta$. In the above derivation we have used the facts: (a) when the eighth power of $\mathcal{Z}$ is developed, if a term $\frac{\eta_x - \rho}{\sigma(\rho)}$ appears with power one, the corresponding expectation is zero; (b) $\sigma^2(\rho) \leq C\rho$ for some $C > 0$ (see Proposition 2.3).

- Suppose $|\mathcal{Z}| \leq \beta$. In this case

$$\left| \frac{\bar{\eta}_{C_1}}{l} - \rho \right| \leq \beta\sigma(\rho).$$

Thus, by the standard estimate on the remainder of a Taylor expansion,

$$(4.21) \qquad |F| \leq \sup\{\alpha''(s) : |s - \rho| \leq \beta\sigma(\rho)\} \mathcal{Z}^2\sigma^2(\rho) \leq C(\beta)\sigma(\rho)\mathcal{Z}^2,$$

where the estimate

$$(4.22) \qquad \sup\{\alpha''(s) : |s - \rho| \leq \beta\sigma(\rho)\} \leq \frac{D(\beta)}{\sigma(\rho)},$$

for every $\rho \geq \rho_0$ and some $D(\beta) > 0$, comes from Lemma 3.5 in [8]. It follows that

$$(4.23) \qquad \mu_\rho[F^4 \mathbf{1}(|\mathcal{Z}| \leq \beta)] \leq C(\beta)\sigma^4(\rho)\mu_\rho[\mathcal{Z}^8] \leq \frac{C(\beta)\rho^2}{l^4},$$

where $\mu_\rho[\mathcal{Z}^8]$ is estimated as in (4.20).

Thus, by (4.20) and (4.23),

$$(4.24) \qquad \mu_\rho[F^4] \leq \frac{C\rho^2}{l^4},$$

for some $C > 0$.

We now estimate the other term in (4.18), namely, $\mu_\rho[e^{2tl|F|}]$. We use again the rough estimate

$$|F| \leq C\sigma(\rho)|\mathcal{Z}|.$$

We get

$$(4.25) \quad \mu_\rho[e^{2tl|F|}] \leq \mu_\rho[e^{Ctl\sigma(\rho)|\mathcal{Z}|}] \leq \mu_\rho[e^{Ct\sum_{x \in C_1}|\eta_x - \rho|}] \leq \mu_\rho[e^{Ctl|\eta_0 - \rho|}].$$

A simple consequence of Condition 2.2 is that there exists a constant $C > 0$ such that

$$|n - \rho| \leq C(|c(n) - c(\rho)| + 1),$$

where $c(\rho)$, $\rho > 0$, is obtained by linear interpolation from $c(n)$, $n \in \mathbb{N}$. Moreover, Lemma 5.8 in [8] states that

$$|c(\rho) - \alpha(\rho)| \leq C\sqrt{\rho}$$

for some $C > 0$ and every $\rho > 0$. Thus,

$$(4.26) \qquad |n - \rho| \leq C|c(n) - \alpha(\rho)| + C\sqrt{\rho} + C.$$

By (4.25), (4.26), the inequality $e^{|x|} \leq e^x + e^{-x}$ and Lemma 4.1, we get

$$(4.27) \qquad \mu_\rho[e^{2tl|F|}] \leq e^{Ctl\sqrt{\rho}}\mu_\rho[e^{Ctl|c(\eta_0) - \alpha(\rho)|}] \leq 2e^{Ctl\sqrt{\rho}}e^{Ct^2l^2\rho e^{Ctl}}.$$

Now, we insert all our estimates in (4.18), and we obtain (4.16).

### 4.1.3. *Final form of the entropy inequality.*  Note that the estimate given in (4.16) for the second factor of (4.10) is worse than the one given in (4.11) for the first factor of (4.10). We can therefore go back to (4.10) and obtain

$$(4.28) \qquad \begin{aligned} &\mu_\rho\left[\exp\left(t\left(\nu_{C_1}^{\bar{\eta}_{C_1}}\left[\sum_{x \in C_1}\tilde{c}(\eta_x)\right] - \mu_\rho\left[\sum_{x \in C_1}\tilde{c}(\eta_x)\right]\right)\right)\right] \\ &\leq \exp[Ct^2\rho e^{Ctl}\sqrt{\rho}e^{Ct^2l^2\rho e^{Ctl}}]. \end{aligned}$$

Now, using (4.5), (4.7) and (4.9), we obtain

$$(4.29) \qquad \begin{aligned} &\nu_\Lambda^N\left[\exp\left(t\sum_{k=1}^m\left(\nu_{C_k}^{\bar{\eta}_{C_k}}\left[\sum_{x \in C_k}\tilde{c}(\eta_x)\right] - \nu_\Lambda^N\left[\sum_{x \in C_k}\tilde{c}(\eta_x)\right]\right)\right)\right] \\ &\leq Ce^{Ct\sqrt{N}}\exp\left[Ct^2\frac{N}{l}e^{Ctl}\sqrt{\rho}e^{Ct^2l^2\rho e^{Ctl}}\right]. \end{aligned}$$

With this estimate, the entropy inequality (4.4) becomes

$$(4.30) \qquad \begin{aligned} &\nu_\Lambda^N\left[f, \sum_{k=1}^m\nu_{C_k}^{\bar{\eta}_{C_k}}\left[\sum_{x \in C_k}\tilde{c}(\eta_x)\right]\right] \\ &\leq \frac{C}{t} + C\sqrt{N} + Ct\frac{N}{l}e^{Ctl}\sqrt{\rho}e^{Ct^2l^2\rho e^{Ctl}} + \frac{1}{t}\,\mathrm{Ent}_{\nu_\Lambda^N}(f). \end{aligned}$$

We now note that the whole argument leading to (4.30) is insensitive to the replacement of $\sum_{k=1}^m \nu_{C_k}^{\bar{\eta}_{C_k}}[\sum_{x \in C_k}\tilde{c}(\eta_x)]$ with its additive inverse. Inequality (4.30) holds therefore for the absolute value of the left-hand side, which gives, for some possibly different $C > 0$,

$$
\text{(4.31)} \quad
\begin{aligned}
&\nu_\Lambda^N \left[ f, \sum_{k=1}^m \nu_{C_k}^{\bar\eta_{C_k}} \left[ \sum_{x \in C_k} \tilde c(\eta_x) \right] \right]^2 \\
&\qquad \leq \frac{C}{t^2} + CN + Ct^2 \frac{N^2}{l^2} e^{Ctl\sqrt\rho} e^{Ct^2 l^2 \rho e^{Ctl}} + \frac{C}{t^2} \operatorname{Ent}^2_{\nu_\Lambda^N}(f).
\end{aligned}
$$

4.1.4. *Optimization of the entropy inequality.* For the optimization of (4.31), we set

$$
t_*^2 := \frac{1 \vee [M \operatorname{Ent}_{\nu_\Lambda^N}(f)]}{N},
$$

where $M$ is a large constant that will be chosen later. We must consider two different regimes:

CASE 1.  $t_* \leq \frac{M}{l\sqrt\rho} \wedge M$.

CASE 2.  $t_* > \frac{M}{l\sqrt\rho} \wedge M$.

Case 1 is more or less modifications of the corresponding argument in the one-block estimate. Case 2 is much more delicate. We use many ideas contained in [2], but we have additional problems due to the unboundedness of particle numbers. The possibility of getting good estimates depends, as we shall see, on the fact that the $\rho$-dependence in (4.31) come through the product $t\sqrt\rho$, and not $t\rho$, that would have required much less work. This justifies the struggle in the previous pages.

CASE 1.  $t_* \leq \frac{M}{l\sqrt\rho} \wedge M$.

Inserting $t_*$ in (4.31), we get, for some possibly different $C > 0$,

$$
\text{(4.32)} \quad
\begin{aligned}
&\nu_\Lambda^N \left[ f, \sum_{k=1}^m \nu_{C_k}^{\bar\eta_{C_k}} \left[ \sum_{x \in C_k} \tilde c(\eta_x) \right] \right]^2 \\
&\qquad \leq CN + C \left( \frac{M}{l} e^{CM} e^{CM^2 e^{CM}} + \frac{1}{M} \right) N \operatorname{Ent}_{\nu_\Lambda^N}(f).
\end{aligned}
$$

By suitably choosing $l$ and $M$, for a given $\varepsilon > 0$, this last term is smaller than

$$
\text{(4.33)} \quad C_\varepsilon N + \varepsilon N \operatorname{Ent}_{\nu_\Lambda^N}(f).
$$

So far we have assumed $f$ normalized, that is, $\nu_\Lambda^N[f] = 1$. For the general case, we may replace $f$ by $f/\nu_\Lambda^N[f]$, and we obtain an inequality stronger that (3.37), for this Case 1.

CASE 2.   $t_* > \frac{M}{l\sqrt{\rho}} \wedge M$.

Suppose, first, $t_* > M$, that is, $1 < \frac{\mathrm{Ent}_{\nu_\Lambda^N}(f)}{MN}$. Recalling $\nu_\Lambda^N(f) = 1$, rough estimates give

$$\nu_\Lambda^N\left[f, \sum_{k=1}^m \nu_{C_k}^{\bar{\eta}_{C_k}}\left[\sum_{x\in C_k} \tilde{c}(\eta_x)\right]\right]^2 \leq CN^2 \leq \frac{CN}{M}\mathrm{Ent}_{\nu_\Lambda^N}(f),$$

that is better than (4.32). We are thus left with the hardest case, $\frac{M}{l\sqrt{\rho}} < t_* \leq M$. We have

$$\left|\nu_\Lambda^N\left[f, \sum_{k=1}^m \nu_{C_k}^{\bar{\eta}_{C_k}}\left[\sum_{x\in C_k} \tilde{c}(\eta_x)\right]\right]\right|$$

$$(4.34) \qquad \leq \left|\nu_\Lambda^N\left[f, \sum_{k=1}^m \mu_{\bar{\eta}_{C_k}/l}\left[\sum_{x\in C_k} \tilde{c}(\eta_x)\right]\right]\right|$$

$$(4.35) \qquad + \left|\nu_\Lambda^N\left[f, \sum_{k=1}^m \left(\nu_{C_k}^{\bar{\eta}_{C_k}}\left[\sum_{x\in C_k} \tilde{c}(\eta_x)\right] - \mu_{\bar{\eta}_{C_k}/l}\left[\sum_{x\in C_k} \tilde{c}(\eta_x)\right]\right)\right]\right|.$$

We begin by estimating (4.35), which is the easiest of the previous two terms. By (4.8), the expression in (4.35) in bounded above by

$$(4.36) \qquad\qquad C\nu_\Lambda^N\left[f\sum_{k=1}^m \sqrt{1 + \frac{\bar{\bar{\eta}}_{C_k}}{l}}\right].$$

By concavity of $\sqrt{\cdot}$,

$$\sum_{k=1}^m \sqrt{1 + \frac{\bar{\bar{\eta}}_{C_k}}{l}} \leq m\sqrt{1 + \frac{N}{L}}.$$

Since $\frac{M}{l\sqrt{\rho}} \leq t_*$ implies $1 \leq l\sqrt{\frac{\rho}{MN}}\sqrt{\mathrm{Ent}_{\nu_\Lambda^N}(f)}$, we get

$$(4.37) \qquad (4.35) \leq Cml\sqrt{1+\rho}\sqrt{\frac{\rho}{MN}}\sqrt{\mathrm{Ent}_{\nu_\Lambda^N}(f)}$$

$$\leq C'L\rho\frac{\sqrt{\mathrm{Ent}_{\nu_\Lambda^N}(f)}}{\sqrt{MN}} \leq \frac{C'}{\sqrt{M}}\sqrt{N\,\mathrm{Ent}_{\nu_\Lambda^N}(f)}.$$

We now turn to the estimate of (4.34). This will require a large number of intermediate steps. We begin by writing

$$\left| \nu_\Lambda^N \left[ f, \sum_{k=1}^m \mu_{\bar\eta_{C_k}/l} \left[ \sum_{x \in C_k} \tilde{c}(\eta_x) \right] \right] \right|$$

$$(4.38) \qquad \leq \left| \nu_\Lambda^N \left[ f \cdot \sum_{k=1}^m \left( \mu_{\bar\eta_{C_k}/l} \left[ \sum_{x \in C_k} \tilde{c}(\eta_x) \right] - \mu_\rho \left[ \sum_{x \in C_k} \tilde{c}(\eta_x) \right] \right) \right] \right|$$

$$(4.39) \qquad + \left| \nu_\Lambda^N \left[ \sum_{k=1}^m \left( \mu_{\bar\eta_{C_k}/l} \left[ \sum_{x \in C_k} \tilde{c}(\eta_x) \right] - \mu_\rho \left[ \sum_{x \in C_k} \tilde{c}(\eta_x) \right] \right) \right] \right|.$$

We now estimate (4.38). From the same argument, the estimate of (4.39) will come for free. Note that

$$(4.40) \qquad \mu_{\bar\eta_{C_k}/l} \left[ \sum_{x \in C_k} \tilde{c}(\eta_x) \right] - \mu_\rho \left[ \sum_{x \in C_k} \tilde{c}(\eta_x) \right] = l F_k,$$

where

$$F_k = \alpha \left( \frac{\bar\eta_{C_k}}{l} \right) - \alpha(\rho) - \alpha'(\rho) \left( \frac{\bar\eta_{C_k}}{l} - \rho \right).$$

We now use essentially the same argument as for the estimate of the right-hand side of (4.17). Set

$$\mathcal{Z}_k := \frac{1}{l} \sum_{x \in C_k} \frac{\eta_x - \rho}{\sigma(\rho)}.$$

We have seen above [see (4.21)] that, for a given $\beta > 0$, the inequality

$$(4.41) \qquad |F_k| \leq C(\beta) \sigma(\rho) |\mathcal{Z}_k|^2$$

holds on the set $\{|\mathcal{Z}_k| \leq \beta\}$, where $C(\beta) \uparrow +\infty$ as $\beta \uparrow +\infty$. For the case $|\mathcal{Z}_k| > \beta$, we write

$$F_k = \alpha(\rho + \sigma(\rho)\mathcal{Z}_k) - \alpha(\rho) - \alpha'(\rho)\sigma(\rho)\mathcal{Z}_k =: H(\mathcal{Z}_k).$$

Since $\alpha'$ is bounded (see formulas (1.3) and (5.2) in [8]), $|H'(z)| \leq C\sigma(\rho)$ for some $C > 0$, so that

$$(4.42) \qquad |F_k| \leq C\sigma(\rho)|\mathcal{Z}_k|.$$

It follows from (4.41) and (4.42) that, if we choose $\beta$ so that $C(\beta)\beta > C$ and define

$$(4.43) \qquad G(z) := \begin{cases} C(\beta)z^2, & \text{for } |z| \leq \beta, \\ \beta^2 C(\beta) + 2\beta C(\beta)(|z| - \beta), & \text{for } |z| > \beta, \end{cases}$$

we have

$$(4.44) \qquad |F_k| \le \sigma(\rho) G(\mathcal{Z}_k).$$

Note that $G(\cdot)$ is a $\mathcal{C}^1$ convex function, and, for some $C > 0$, $G(z) \le C(|z| \wedge z^2)$. By convexity of $G(\cdot)$,

$$\sum_{k=1}^m G(\mathcal{Z}_k) \le \frac{1}{m} \sum_{k,p=1}^m G\left(\frac{\bar\eta_{C_k} - \bar\eta_{C_p}}{l\sigma(\rho)}\right).$$

Thus,

$$(4.38) \le \frac{l}{m}\sigma(\rho) \sum_{k,p=1}^m \nu_\Lambda^N\left[f \cdot G\left(\frac{\bar\eta_{C_k} - \bar\eta_{C_p}}{l\sigma(\rho)}\right)\right]$$

$$(4.45) \qquad = \frac{l}{m}\sigma(\rho) \sum_{k,p=1}^m \nu_\Lambda^N\left[f, G\left(\frac{\bar\eta_{C_k} - \bar\eta_{C_p}}{l\sigma(\rho)}\right)\right]$$

$$(4.46) \qquad + \frac{l}{m}\sigma(\rho) \sum_{k,p=1}^m \nu_\Lambda^N\left[G\left(\frac{\bar\eta_{C_k} - \bar\eta_{C_p}}{l\sigma(\rho)}\right)\right].$$

To estimate (4.45) and (4.46), we introduce the $\sigma$-fields

$$\mathcal{F}_{k,p} = \sigma\{\eta_x : x \notin C_k \cup C_p\}.$$

Note that

$$\nu_\Lambda^N\left[G\left(\frac{\bar\eta_{C_k} - \bar\eta_{C_p}}{l\sigma(\rho)}\right)\Big|\mathcal{F}_{k,p}\right] = \nu_{C_k \cup C_p}^{\bar\eta_{C_k} + \bar\eta_{C_p}}\left[G\left(\frac{\bar\eta_{C_k} - \bar\eta_{C_p}}{l\sigma(\rho)}\right)\right]$$

$$\le \frac{C}{l^2\sigma^2(\rho)} \nu_{C_k \cup C_p}^{\bar\eta_{C_k} + \bar\eta_{C_p}}[(2\bar\eta_{C_k} - k)^2]|_{k=\bar\eta_{C_k} + \bar\eta_{C_p}},$$

where we have used the inequality $G(z) \le Cz^2$. The distribution of $\bar\eta_{C_k}$ under $\nu_{C_k \cup C_p}^{\bar\eta_{C_k} + \bar\eta_{C_p}}$ has been analyzed in Lemmas 8.2 and 8.6 of [4], where we have obtained uniform Gaussian estimates, that imply

$$(4.47) \qquad \nu_\Lambda^N\left[G\left(\frac{\bar\eta_{C_k} - \bar\eta_{C_p}}{l\sigma(\rho)}\right)\Big|\mathcal{F}_{k,p}\right] \le \frac{C(\bar\eta_{C_k} + \bar\eta_{C_p})}{l^2\sigma^2(\rho)}.$$

In particular,

$$(4.48) \qquad \sum_{k,p=1}^m \nu_\Lambda^N\left[G\left(\frac{\bar\eta_{C_k} - \bar\eta_{C_p}}{l\sigma(\rho)}\right)\Big|\mathcal{F}_{k,p}\right] \le \frac{2CmN}{l^2\sigma^2(\rho)}.$$

Thus,

$$(4.49) \qquad (4.46) \le \frac{2CN}{l\sigma(\rho)} \le C'\sqrt{\frac{N}{M}\,\mathrm{Ent}_{\nu_\Lambda^N}(f)},$$

where the inequality $t_* > \frac{M}{l\sqrt{\rho}}$ has been used.

We now estimate (4.45). By the usual formula for the conditional covariance,

$$\frac{l}{m}\sigma(\rho)\sum_{k,p=1}^{m}\nu_{\Lambda}^{N}\Big[f, G\Big(\frac{\bar{\eta}_{C_k}-\bar{\eta}_{C_p}}{l\sigma(\rho)}\Big)\Big]$$

$$(4.50) \qquad = \frac{l}{m}\sigma(\rho)\sum_{k,p=1}^{m}\nu_{\Lambda}^{N}\Big[f, \nu_{\Lambda}^{N}\Big[G\Big(\frac{\bar{\eta}_{C_k}-\bar{\eta}_{C_p}}{l\sigma(\rho)}\Big)\big|\mathcal{F}_{k,p}\Big]\Big]$$

$$(4.51) \qquad + \frac{l}{m}\sigma(\rho)\sum_{k,p=1}^{m}\nu_{\Lambda}^{N}\Big[\nu_{\Lambda}^{N}\Big[f, G\Big(\frac{\bar{\eta}_{C_k}-\bar{\eta}_{C_p}}{l\sigma(\rho)}\Big)\big|\mathcal{F}_{k,p}\Big]\Big].$$

Using (4.48) and $\nu_{\Lambda}^{N}(f)=1$, we get for (4.50) the same upper bound as for (4.46). The key is now to estimate (4.51). In what follows, for compactness of notation, we write $\nu_{k,p}[\cdot]$ for $\nu_{\Lambda}^{N}[\cdot|\mathcal{F}_{k,p}]$ and $G_{k,p}$ for $G(\frac{\bar{\eta}_{C_k}-\bar{\eta}_{C_p}}{l\sigma(\rho)})$. The covariance

$$\nu_{k,p}[f, G_{k,p}]$$

is left unchanged if we replace $f$ by $f - \nu_{k,p}[\sqrt{f}]^2 = (\sqrt{f}-\nu_{k,p}[\sqrt{f}])(\sqrt{f}+\nu_{k,p}[\sqrt{f}])$. Thus, by the Cauchy–Schwarz inequality,

$$(4.52) \qquad \begin{aligned} |\nu_{k,p}[f, G_{k,p}]| \\ \leq \{\nu_{k,p}[(\sqrt{f}-\nu_{k,p}[\sqrt{f}])^2|G_{k,p}-\nu_{k,p}[G_{k,p}]|] \\ \times \nu_{k,p}[(\sqrt{f}+\nu_{k,p}[\sqrt{f}])^2|G_{k,p}-\nu_{k,p}[G_{k,p}]|]\}^{1/2} \\ =: \{D_1 D_2\}^{1/2}. \end{aligned}$$

In the above expression, $G_{k,p}$ is bounded from above by $C\frac{|\bar{\eta}_{C_k}-\bar{\eta}_{C_p}|}{l\sigma(\rho)}$, for some $C>0$. Its $\nu_{k,p}$-expectation is estimated again via the Gaussian estimates in Lemmas 8.2 and 8.6 of [4], giving

$$\nu_{k,p}[G_{k,p}] \leq C\frac{\sqrt{\bar{\eta}_{C_k}+\bar{\eta}_{C_p}}}{l\sigma(\rho)}.$$

Thus,

$$(4.53) \qquad \begin{aligned} D_1 &\leq C\nu_{k,p}\Big[(\sqrt{f}-\nu_{k,p}[\sqrt{f}])^2\Big(\frac{|\bar{\eta}_{C_k}-\bar{\eta}_{C_p}|}{l\sigma(\rho)}+\frac{\sqrt{\bar{\eta}_{C_k}+\bar{\eta}_{C_p}}}{l\sigma(\rho)}\Big)\Big] \\ &\leq \frac{C}{l\sigma(\rho)}\nu_{k,p}[(\sqrt{f}-\nu_{k,p}[\sqrt{f}])^2] \\ &\quad \times \Big(\frac{1}{t}\log\nu_{k,p}[e^{t|\bar{\eta}_{C_k}-\bar{\eta}_{C_p}|}]+\sqrt{\bar{\eta}_{C_k}+\bar{\eta}_{C_p}}\Big) \end{aligned}$$

$$+ \frac{1}{t} \frac{C}{l\sigma(\rho)} \, \mathrm{Ent}_{\nu_{k,p}}((\sqrt{f} - \nu_{k,p}[\sqrt{f}])^2),$$

where we used the by now usual entropy inequality (3.19). The estimates on the moment generating function

$$\nu_{k,p}[e^{t|\bar{\eta}_{C_k} - \bar{\eta}_{C_p}|}] = \nu_{C_k \cup C_p}^{\bar{\eta}_{C_k} + \bar{\eta}_{C_p}}[e^{2t|\bar{\eta}_{C_k} - (\bar{\eta}_{C_k} + \bar{\eta}_{C_p})/2|}]$$

can be done through arguments that we used already. We therefore only sketch the steps needed:

1. By comparison between ensembles, using (4.6) as in (4.7), we reduce the problem to estimating a moment generating function with respect to the grand canonical measure $\mu_{\rho_{k,p}}$, with $\rho_{k,p} = \frac{\bar{\eta}_{C_k} + \bar{\eta}_{C_p}}{2l}$.

2. We are then left to estimate quantities of the form

$$(4.54) \qquad\qquad \mu_{\rho_{k,p}}[e^{t|\eta_x - \rho_{k,p}|}].$$

3. Lemma 4.1 allows us to estimate

$$(4.55) \qquad\qquad \mu_{\rho_{k,p}}[e^{t|c(\eta_x) - \alpha(\rho_{k,p})|}] \le 2e^{C\rho_{k,p}t^2 e^{Ct}},$$

   for some $C > 0$. Now, for $r > 0$, let $c(r)$ be obtained by linear interpolation of $c(n)$, $n \in \mathbb{N}$. A simple consequence of Condition 2.2 is that there exists a constant $C > 0$ such that

$$(4.56) \qquad\qquad |r - s| \le C(|c(r) - c(s)| + 1)$$

   for every $r, s > 0$. Moreover, by Lemma 5.8 in [8],

$$(4.57) \qquad\qquad |c(r) - \alpha(r)| \le C\sqrt{r}$$

   for every $r > 0$. Putting (4.56) and (4.57) together, we get

$$(4.58) \qquad\qquad |\eta_x - \rho_{k,p}| \le C|c(\eta_x) - \alpha(\rho_{k,p})| + C\sqrt{\bar{\eta}_{C_k} + \bar{\eta}_{C_p}}.$$

4. Using (4.55) and (4.58), we get

$$(4.59) \quad \nu_{k,p}[e^{t|\bar{\eta}_{C_k} - \bar{\eta}_{C_p}|}] \le C \exp[Ct\sqrt{\bar{\eta}_{C_k} + \bar{\eta}_{C_p}} + Ct^2(\bar{\eta}_{C_k} + \bar{\eta}_{C_p})e^{Ct}].$$

Collecting these last estimates and (4.53), we obtain, for $t \le 1$,

$$(4.60) \quad \begin{aligned} D_1 &\le \frac{C}{l\sigma(\rho)} \nu_{k,p}[\sqrt{f}, \sqrt{f}] \left[ \frac{1}{t} + \sqrt{\bar{\eta}_{C_k} + \bar{\eta}_{C_p}} + t(\bar{\eta}_{C_k} + \bar{\eta}_{C_p}) \right] \\ &\quad + \frac{1}{t} \frac{C}{l\sigma(\rho)} \mathrm{Ent}_{\nu_{k,p}}((\sqrt{f} - \nu_{k,p}[\sqrt{f}])^2). \end{aligned}$$

To optimize (4.60), we set

$$t_*^2 := \frac{1 \vee (\mathrm{Ent}_{\nu_{k,p}}((\sqrt{f} - \nu_{k,p}[\sqrt{f}])^2))/\nu_{k,p}[\sqrt{f}, \sqrt{f}]}{\bar{\eta}_{C_k} + \bar{\eta}_{C_p}}.$$

If $t_* \leq 1$, inserting $t_*$ in (4.60), we get

$$
(4.61) \quad
\begin{aligned}
D_1 &\leq \frac{C}{l\sigma(\rho)} \nu_{k,p}[\sqrt{f}, \sqrt{f}] \\
&\times \left( \sqrt{\bar{\eta}_{C_k} + \bar{\eta}_{C_p}} + \sqrt{\bar{\eta}_{C_k} + \bar{\eta}_{C_p}} \sqrt{\frac{\mathrm{Ent}_{\nu_{k,p}}((\sqrt{f} - \nu_{k,p}[\sqrt{f}])^2)}{\nu_{k,p}[\sqrt{f}, \sqrt{f}]}} \right).
\end{aligned}
$$

On the other hand, for $t_* > 1$,

$$
\begin{aligned}
D_1 &\leq \frac{C(\bar{\eta}_{C_k} + \bar{\eta}_{C_p})}{l\sigma(\rho)} \nu_{k,p}[\sqrt{f}, \sqrt{f}] \\
&\leq \frac{C\sqrt{\bar{\eta}_{C_k} + \bar{\eta}_{C_p}}}{l\sigma(\rho)} \sqrt{\nu_{k,p}[\sqrt{f}, \sqrt{f}] \, \mathrm{Ent}_{\nu_{k,p}}((\sqrt{f} - \nu_{k,p}[\sqrt{f}])^2)},
\end{aligned}
$$

that implies (4.61). Now, let $\mathcal{E}_{k,p}(\cdot, \cdot)$ be the Dirichlet form of the zero-range process in $C_k \cup C_p$, where one point of $C_k$ is assumed to be nearest neighbor to one of $C_p$, so that exchange of particles is allowed. By (3.18), this process satisfies a logarithmic Sobolev inequality with a constant $s(l)$ depending only on $l$. Therefore,

$$
(4.62) \quad \mathrm{Ent}_{\nu_{k,p}}((\sqrt{f} - \nu_{k,p}[\sqrt{f}])^2) \leq s(l)\mathcal{E}_{k,p}(\sqrt{f}, \sqrt{f}).
$$

Thus, by (4.61) and (4.62),

$$
(4.63) \quad
\begin{aligned}
D_1 &\leq \frac{C\sqrt{\bar{\eta}_{C_k} + \bar{\eta}_{C_p}}}{l\sigma(\rho)} \nu_{k,p}[\sqrt{f}, \sqrt{f}] \\
&+ \frac{C\sqrt{\bar{\eta}_{C_k} + \bar{\eta}_{C_p}}}{l\sigma(\rho)} \sqrt{\nu_{k,p}[\sqrt{f}, \sqrt{f}] s(l)\mathcal{E}_{k,p}(\sqrt{f}, \sqrt{f})}.
\end{aligned}
$$

Going back to (4.52), we have to estimate $D_2$. This estimate goes along the same lines as the one for $D_1$, giving

$$
(4.64) \quad
\begin{aligned}
D_2 &\leq \frac{C\sqrt{\bar{\eta}_{C_k} + \bar{\eta}_{C_p}}}{l\sigma(\rho)} \nu_{k,p}[f] \\
&+ \frac{C\sqrt{\bar{\eta}_{C_k} + \bar{\eta}_{C_p}}}{l\sigma(\rho)} \sqrt{\nu_{k,p}[f] s(l)\mathcal{E}_{k,p}(\sqrt{f}, \sqrt{f})}.
\end{aligned}
$$

Now we consider the product $D_1 D_2$. We can use the spectral gap estimate $\nu_{k,p}[\sqrt{f}, \sqrt{f}] \leq s(l)\mathcal{E}_{k,p}(\sqrt{f}, \sqrt{f})$ and the trivial bound $\nu_{k,p}[\sqrt{f}, \sqrt{f}] \leq \nu_{k,p}[f]$ to obtain, for some new constant $C$ that may depend on $l$,

$$
(4.65) \quad D_1 D_2 \leq \frac{C(\bar{\eta}_{C_k} + \bar{\eta}_{C_p})}{\sigma^2(\rho)} \nu_{k,p}[f] \mathcal{E}_{k,p}(\sqrt{f}, \sqrt{f}).
$$

So, by the Cauchy–Schwarz inequality and the fact that $\bar{\eta}_{C_k} + \bar{\eta}_{C_p}$ is $\mathcal{F}_{k,p}$-measurable, we get

$$\nu_\Lambda^N[\sqrt{D_1 D_2}] \leq \frac{C}{\sigma(\rho)} \sqrt{\nu_\Lambda^N[(\bar{\eta}_{C_k} + \bar{\eta}_{C_p})\nu_{k,p}[f]]\nu_\Lambda^N[\mathcal{E}_{k,p}(\sqrt{f},\sqrt{f})]}$$

(4.66)

$$= \frac{C}{\sigma(\rho)} \sqrt{\nu_\Lambda^N[(\bar{\eta}_{C_k} + \bar{\eta}_{C_p})f]\nu_\Lambda^N[\mathcal{E}_{k,p}(\sqrt{f},\sqrt{f})]}.$$

Thus, if $C$ is some constant possibly dependent on $l$ which may change from step to step, by (4.52), (4.66) and Jensen's inequality, we have

$$|(4.51)| \leq \frac{l}{m}\sigma(\rho) \sum_{k,p=1}^m \nu_\Lambda^N[\sqrt{D_1 D_2}]$$

(4.67)

$$\leq C\frac{1}{m} \sum_{k,p=1}^m \sqrt{\nu_\Lambda^N[(\bar{\eta}_{C_k} + \bar{\eta}_{C_p})f]\nu_\Lambda^N[\mathcal{E}_{k,p}(\sqrt{f},\sqrt{f})]}$$

$$\leq C\frac{1}{m} \sum_{k,p=1}^m \sqrt{\nu_\Lambda^N[\bar{\eta}_{C_k}f]\nu_\Lambda^N[\mathcal{E}_{k,p}(\sqrt{f},\sqrt{f})]}$$

$$\leq C \sum_k \sqrt{\nu_\Lambda^N[\bar{\eta}_{C_k}f]\frac{1}{m}\sum_p \nu_\Lambda^N[\mathcal{E}_{k,p}(\sqrt{f},\sqrt{f})]}.$$

Now consider the term

$$\nu_\Lambda^N[\mathcal{E}_{k,p}(\sqrt{f},\sqrt{f})] = \sum_{x \in C_k} \sum_{y \in C_k, y \sim x} \nu_\Lambda^N[c(\eta_x)(\partial_{xy}\sqrt{f})^2]$$

(4.68)

$$+ \sum_{x \in C_p} \sum_{y \in C_p, y \sim x} \nu_\Lambda^N[c(\eta_x)(\partial_{xy}\sqrt{f})^2]$$

$$+ \nu_\Lambda^N[(c(\eta_{x_k}) + c(\eta_{x_p}))(\partial_{x_k,x_p}\sqrt{f})^2],$$

where $x_k$ and $x_p$ are the two sites in, respectively, $C_k$ and $C_p$ that may exchange particles. The first two summand in (4.68) are both bounded above by $\mathcal{E}_{\nu_\Lambda^N}(\sqrt{f},\sqrt{f})$. To bound the third summand, we use the following standard argument. Assume $x_p = x_k + h$, $h > 0$. Then, by Jensen's inequality,

$$\nu_\Lambda^N[c(\eta_{x_k})(\partial_{x_k,x_p}\sqrt{f})^2]$$

(4.69)

$$= \nu_\Lambda^N\left[c(\eta_{x_k})\left(\sum_{i=1}^h \partial_{x_k+i-1,x_k+i}\sqrt{f}(\eta - \delta_{x_k} + \delta_{x_k+i-1})\right)^2\right]$$

$$\leq h\sum_{i=1}^h \nu_\Lambda^N[c(\eta_{x_k})(\partial_{x_k+i-1,x_k+i}\sqrt{f}(\eta - \delta_{x_k} + \delta_{x_k+i-1}))^2]$$

$$= h \sum_{i=1}^{h} \nu_\Lambda^N [c(\eta_{x_k+i-1})(\partial_{x_k+i-1,x_k+i} \sqrt{f}(\eta))^2],$$

where the last equality is a simple consequence of the detail balance (2.3). Clearly, since $h \leq L$, this last term is bounded above by $L\mathcal{E}_{\nu_\Lambda^N}(\sqrt{f}, \sqrt{f})$. Putting it all together, we have

$$\sum_p \nu_\Lambda^N [\mathcal{E}_{k,p}(\sqrt{f}, \sqrt{f})] \leq CmL\mathcal{E}_{\nu_\Lambda^N}(\sqrt{f}, \sqrt{f}) \tag{4.70}$$

for some $C > 0$. Using (4.70), we can continue (4.67), obtaining, using Jensen's inequality again and $\nu_\Lambda^N(f) = 1$,

$$\begin{aligned}
|(4.51)| &\leq C \sum_k \sqrt{\nu_\Lambda^N[\bar\eta_{C_k} f] L \mathcal{E}_{\nu_\Lambda^N}(\sqrt{f}, \sqrt{f})} \\
&\leq Cm \sqrt{\nu_\Lambda^N \left[\frac{1}{m} \sum_k \bar\eta_{C_k} f\right] L \mathcal{E}_{\nu_\Lambda^N}(\sqrt{f}, \sqrt{f})} \\
&\leq C\sqrt{N} L \sqrt{\mathcal{E}_{\nu_\Lambda^N}(\sqrt{f}, \sqrt{f})}.
\end{aligned} \tag{4.71}$$

Now, going back to (4.45), (4.46), (4.49), (4.50), (4.51) and (4.71), we have

$$(4.38) \leq \sqrt{N} \left[\frac{C\sqrt{\mathrm{Ent}_{\nu_\Lambda^N}(f)}}{\sqrt{M}} + D_l L \sqrt{\mathcal{E}_{\nu_\Lambda^N}(\sqrt{f}, \sqrt{f})}\right], \tag{4.72}$$

where $D_l$ depends on $l$, but $C$ does not. This is very important since, in this whole argument, $l$ and $M$ are taken large, but not in an independent way. This last estimate easily extends to (4.39): in this case only the term (4.46) survives. Finally, putting together (4.37) and (4.72), and choosing $M$ large enough, we have

$$\begin{aligned}
\nu_\Lambda^N &\left[f, \sum_{k=1}^{m} \nu_{C_k}^{\bar\eta_{C_k}} \left[\sum_{x \in C_k} \tilde{c}(\eta_x)\right]\right]^2 \\
&\leq N[C_\varepsilon L^2 \mathcal{E}_{\nu_\Lambda^N}(\sqrt{f}, \sqrt{f}) + \varepsilon \, \mathrm{Ent}_{\nu_\Lambda^N}(f)].
\end{aligned} \tag{4.73}$$

This completes the analysis of Case 2. Together with the estimates obtained for Case 1, the proof of Proposition 3.10 is completed for densities away from zero.

### 4.2. *Case of small density.*

For the case of small density, the two blocks argument is irrelevant. We indeed succeed in proving Corollary 3.11 directly. We restate Corollary 3.11 in more convenient terms for small densities.

Lemma 4.4. *For every $\varepsilon > 0$, there exists $\rho_0 > 0$ and a constant $C_\varepsilon > 0$ such that, for $\frac{N}{L} \le \rho_0$,*

$$(4.74) \qquad \nu_\Lambda^N\left[f, \sum_{x \in \Lambda} c(\eta_x)\right]^2 \le N\nu_\Lambda^N[f](C_\varepsilon \nu_\Lambda^N[f] + \varepsilon \operatorname{Ent}_{\nu_\Lambda^N}(f)).$$

Proof. Without loss of generality, we again assume $\nu_\Lambda^N[f] = 1$. One more time we begin with the entropy inequality. Letting $\rho = N/L$, we have

$$
\begin{aligned}
(4.75) \qquad &\nu_\Lambda^N\left[f, \sum_{x \in \Lambda} c(\eta_x)\right] \\
&= \nu_\Lambda^N\left[f, \sum_{x \in \Lambda}\left(c(\eta_x) - \frac{\alpha(\rho)}{\rho}\eta_x\right)\right] \\
&\le \frac{1}{t}\log \nu_\Lambda^N\left[\exp\left(t\sum_{x \in \Lambda}\left(c(\eta_x) - \frac{\alpha(\rho)}{\rho}\eta_x - \nu_\Lambda^N\left[c(\eta_x) - \frac{\alpha(\rho)}{\rho}\eta_x\right]\right)\right)\right] \\
&\quad + \frac{1}{t}\operatorname{Ent}_{\nu_\Lambda^N}(f).
\end{aligned}
$$

The next step consists in using the comparison between ensembles as in (4.5)–(4.13): we do not repeat the argument. Taking into account that

$$\mu_\rho\left[c(\eta_x) - \frac{\alpha(\rho)}{\rho}\eta_x\right] = 0,$$

we are lead to the following estimate:

$$
\begin{aligned}
(4.76) \qquad &\nu_\Lambda^N\left[\exp\left(t\sum_{x \in \Lambda}\left(c(\eta_x) - \frac{\alpha(\rho)}{\rho}\eta_x - \nu_\Lambda^N\left[c(\eta_x) - \frac{\alpha(\rho)}{\rho}\eta_x\right]\right)\right)\right] \\
&\le Ce^{C\sqrt{N}t}(\mu_\rho[e^{t(c(\eta_x) - (\alpha(\rho)/\rho)\eta_x)}])^L.
\end{aligned}
$$

The point is now to estimate the moment generating function

$$\mu_\rho[e^{t(c(\eta_x) - (\alpha(\rho)/\rho)\eta_x)}].$$

By an extension of Lemma 4.1, we would get an estimate of the form $\exp[C\alpha t^2 e^{c|t|}]$, where $\alpha = \alpha(\rho)$. Recalling that $\alpha$ and $\rho$ have the same order as $\rho \downarrow 0$, this estimate is not good for $\rho$ (or $\alpha$) small. We now show that the right order, for small $\alpha$ and $t$, is $e^{C\alpha^2 t^2}$, and not $e^{C\alpha t^2}$. Indeed, let

$$F(t, \alpha) := \mu_\rho[e^{t(c(\eta_x) - (\alpha/\rho(\alpha))\eta_x)}],$$

where we choose to take $\alpha$ rather than $\rho = \rho(\alpha)$ as main parameter. $F(t, \alpha)$ is an entire analytic function, since $\frac{\alpha}{\rho(\alpha)}$ extends analytically at $\alpha = 0$ where it takes the value $c(1)$, as shown by a direct computation. Note that $F$ is of the form

$$(4.77) \qquad F(t, \alpha) = \mu_\rho[e^{t\varphi(\alpha,n)}],$$

with

$$\varphi(\alpha, n) = c(n) - \frac{\alpha}{\rho(\alpha)} n.$$

We have

$$(4.78) \qquad F(0, \alpha) \equiv 1,$$

$$(4.79) \qquad \partial_t F(0, \alpha) = \mu_\rho(\varphi(\alpha, n)) \equiv 0,$$

so that

$$(4.80) \qquad \partial_t \partial_\alpha^k F(0, 0) = 0$$

for all $k \geq 0$. Finally,

$$\partial_t^2 F(0, \alpha) = \mu_\rho[\varphi^2(\alpha, n)],$$

from which it follows that

$$
\begin{aligned}
(4.81) \qquad \partial_\alpha \partial_t^2 F(0, \alpha) = &-\frac{Z'(\alpha)}{Z(\alpha)} \mu_\rho[\varphi^2(\alpha, n)] \\
&+ \frac{1}{Z(\alpha)} \sum_{n=1}^{+\infty} n\varphi^2(\alpha, n) \frac{\alpha^{n-1}}{c(n)!} \\
&+ 2\mu_\rho[\varphi(\alpha, n) \, \partial_\alpha \varphi(\alpha, n)].
\end{aligned}
$$

It is therefore easy to see that

$$(4.82) \qquad \partial_\alpha \partial_t^2 F(0, 0) = 0,$$

provided

$$(4.83) \qquad \varphi(0, 0) = \varphi(0, 1) = 0,$$

which is seen by direct inspection.

Fix now a constant $M > 0$. By (4.78), (4.80) and (4.82), it follows that the inequality

$$(4.84) \qquad F(t, \alpha) \leq e^{C_M \alpha^2 t^2}$$

holds for every $t \in [0, M]$ and $\rho \leq 1$, where $C_M$ is a constant that may depend on $M$. Thus, using (4.75), (4.76) and (4.84), we have

$$(4.85) \qquad \nu_\Lambda^N \left[ f, \sum_{x \in \Lambda} c(\eta_x) \right] \leq \frac{C}{t} + C\sqrt{N} + C_M N \rho t + \frac{1}{t} \operatorname{Ent}_{\nu_\Lambda^N}(f),$$

for $t \in [0, M]$ and $\rho \leq 1$. We again observe that the same estimate is obtained if $\sum_{x \in \Lambda} c(\eta_x)$ is replaced by its opposite. It follows that, for some possibly different $C, C_M$,

$$(4.86) \qquad \nu_\Lambda^N \left[ f, \sum_{x \in \Lambda} c(\eta_x) \right]^2 \leq \frac{C}{t^2} + CN + C_M N^2 \rho^2 t^2 + \frac{1}{t^2} \operatorname{Ent}_{\nu_\Lambda^N}^2(f).$$

For optimizing (4.86), we set

$$(4.87) \qquad t_*^2 = \frac{1 \vee M[\operatorname{Ent}_{\nu_\Lambda^N}(f)]}{N}.$$

If $t_* \leq M$, then we can plug it in (4.86), obtaining

$$(4.88) \qquad \begin{aligned} &\nu_\Lambda^N \left[ f, \sum_{x \in \Lambda} c(\eta_x) \right]^2 \\ &\qquad \leq 2CN + C_M N \rho^2 (1 \vee [M \operatorname{Ent}_{\nu_\Lambda^N}(f)]) + \frac{N}{M} \operatorname{Ent}_{\nu_\Lambda^N}(f). \end{aligned}$$

In the case $t_* > M$, we have

$$1 \leq \frac{1}{MN} \operatorname{Ent}_{\nu_\Lambda^N}(f),$$

and so

$$(4.89) \qquad \nu_\Lambda^N \left[ f, \sum_{x \in \Lambda} c(\eta_x) \right]^2 \leq CN^2 \leq \frac{CN}{M} \operatorname{Ent}_{\nu_\Lambda^N}(f).$$

Using (4.88) and (4.89), choosing first $M$ large enough and then $\rho_0$ small enough, the statement of Lemma 4.4 follows, for $\nu_\Lambda^N[f] = 1$. The general case comes easily applying this last result to $f/\nu_\Lambda^N[f]$. $\quad \square$

## 5. Two blocks estimates II: proof of Proposition 3.12.

We rewrite (3.38) in the form

$$(5.1) \qquad \begin{aligned} &\nu_\Lambda^N \left[ f, \sum_{x \in \Lambda} \rho h(\eta_x) \right]^2 \\ &\qquad \leq N \nu_\Lambda^N[f][C_\varepsilon \nu_\Lambda^N[f] + C_\varepsilon L^2 \mathcal{E}_{\nu_\Lambda^N}(\sqrt{f}, \sqrt{f}) + \varepsilon \operatorname{Ent}_{\nu_\Lambda^N}(f)]. \end{aligned}$$

This is the same as (3.37), with $\rho h(\eta_x)$ in place of $c(\eta_x)$. We now simply follow the whole proof of (3.37), and show that the replacement of $c(\eta_x)$ with $\rho h(\eta_x)$ does not cause any harm, up to minor modifications.

The first part of the argument (3.32)–(3.35) goes through with no changes, except that, in proving (3.35) with $\rho h(\eta_x)$ in place of $c(\eta_x)$, we use the estimates in (3.22) rather than the ones in (3.21). Thus, the statement that corresponds to Corollary 3.9 follows.

Thus, all we have to show is that there is choice of $\delta = \delta(\rho)$ such that the statement of Proposition 3.10 holds replacing $\tilde{c}(\eta_x)$ with $\rho h(\eta_x) - \delta \eta_x$.

5.1. *Case of density away from zero.* We now follow the proof given in Section 4.1, for $\rho \geq \rho_0$, for some $\rho_0 > 0$. The steps described in equations (4.4)–(4.7) require no changes, since they are only based on the fact that $\tilde{c}(\eta_x)$ depends only on $\eta_x$. The first step that needs justification is (4.8), which now reads

$$(5.2) \qquad |\nu_\Lambda^N[\rho h(\eta_x)] - \mu_\rho[\rho h(\eta_x)]| \leq \frac{C}{L}\sqrt{1+\rho}.$$

This follows from Corollary 6.1 in [8], provided we prove that

$$(5.3) \qquad \mu_\rho[h(\eta_x), h(\eta_x)] \leq \frac{C}{\rho}$$

for some $C > 0$. Inequality (5.3) can be proved in various ways. Perhaps the fastest one makes use of the Poincaré inequality (1.6)

$$(5.4) \qquad \mu_\rho[f, f] \leq B\mu_\rho[c(n)[f(n-1) - f(n)]^2],$$

that holds for each $f$ having finite $\mu_\rho$-variance, where $C$ is a constant that does not depend on $\rho$. Thus,

$$
\begin{aligned}
(5.5) \qquad \mu_\rho[h(\eta_x), h(\eta_x)] &\leq C\mu_\rho\left[c(n)\left[\frac{n}{c(n)} - \frac{n+1}{c(n+1)}\right]^2\right] \\
&\leq C\mu_\rho\left[c(n)\left[\frac{1}{c(n)c(n+1)}[n(c(n+1) - c(n))]\right]^2\right] \\
&\leq D\mu_\rho\left[\frac{1}{c(n+1)}\right],
\end{aligned}
$$

where, in the last step, we have used the fact that $c(\cdot)$ is Lipschitz and that $\frac{n}{c(n)}$ is bounded and bounded away from zero (Proposition 2.3). Finally, a direct computation shows that

$$(5.6) \qquad \mu_\rho\left[\frac{1}{c(n+1)}\right] = \frac{1}{\alpha}\frac{Z(\alpha) - 1}{Z(\alpha)} \leq \frac{1}{\alpha},$$

which, together with (5.5) and (5.4), proves (5.3). At this point, inequality (4.10) for the new $\tilde{c}$ follows immediately.

The estimate of the first factor of (4.10) simply uses (4.8), so now we use (5.2) instead. The second factor requires some work. As in the estimate of the right-hand side of (4.17), it can be rewritten in the form $\mu_\rho[e^{ltG}]$, where

$$(5.7) \qquad G = \rho\left[\frac{\bar{\eta}_{C_1}/l}{\alpha(\bar{\eta}_{C_1}/l)} - \frac{\rho}{\alpha(\rho)}\right] - \delta\left(\frac{\bar{\eta}_{C_1}}{l} - \rho\right).$$

Set $\gamma(\rho) = \frac{\rho}{\alpha(\rho)}$ and

$$\delta = \rho\gamma'(\rho) = \rho\left[\frac{1}{\alpha(\rho)} - \frac{\rho\alpha'(\rho)}{\alpha^2(\rho)}\right]. \tag{5.8}$$

Now, after replacing $F$ with $G$ in (4.18), we distinguish the two cases $|\mathcal{Z}| > \beta$ and $|\mathcal{Z}| \leq \beta$, for some $\beta > 0$ suitably chosen, where $\mathcal{Z}$ was defined in (4.19).

- For the case $|\mathcal{Z}| > \beta$, in order to obtain the inequality that corresponds to (4.20), we must show

$$|G| \leq C\sigma(\rho)|\mathcal{Z}|. \tag{5.9}$$

We examine separately the two summands in (5.7). For the term $\delta(\bar{\eta}_{C_1}/l - \rho)$, the bound in (5.9) is obvious. For the term

$$\rho\left[\frac{\bar{\eta}_{C_1}/l}{\alpha(\bar{\eta}_{C_1}/l)} - \frac{\rho}{\alpha(\rho)}\right],$$

we set $x = \frac{\bar{\eta}_{C_1}}{l}$. We must show that

$$|\gamma(x) - \gamma(\rho)| \leq \frac{C}{\rho}|x - \rho|. \tag{5.10}$$

Since $\gamma$ is bounded, inequality (5.10) is obvious for $|x - \rho| > \frac{\rho}{2}$. For the case $|x - \rho| \leq \frac{\rho}{2}$, it is enough to observe that

$$\sup\left\{|\gamma'(x)| : x \geq \frac{\rho}{2}\right\} \leq \frac{C}{\rho},$$

as it follows from $\gamma'(\rho) = \frac{\alpha(\rho) - \rho\alpha'(\rho)}{\alpha^2(\rho)}$ and the fact that $\rho$ and $\alpha(\rho)$ have the same order at both zero and $+\infty$ (Proposition 2.3).

- For the case $|\mathcal{Z}| \leq \beta$, we have, similarly to (4.21),

$$|G| \leq \rho\sup\{\gamma''(s) : |s - \rho| \leq \beta\sigma(\rho)\}\mathcal{Z}^2\sigma^2(\rho). \tag{5.11}$$

Note that

$$\gamma''(s) = -\frac{2\alpha'(s)}{\alpha^2(s)} + \frac{2s(\alpha'(s))^2}{\alpha^3(s)} - \frac{s\alpha''(s)}{\alpha^2(s)}. \tag{5.12}$$

Since $\rho \geq \rho_0$, we may assume that $\beta$ is small enough (depending on $\rho_0$) so that $|s - \rho| \leq \beta\sigma(\rho) \Rightarrow s \geq \frac{\rho}{2}$. The first two summands in (5.12) are bounded by $\frac{C}{s^2}$. The third is the dominant one. By (4.22), it is bounded by $\frac{C}{s\sigma(s)}$. Therefore,

$$|G| \leq C(\beta)\sigma(\rho)\mathcal{Z}^2, \tag{5.13}$$

which completes the inequality that corresponds to (4.21).

At this point, the part of the proof from (4.23) to (4.31) requires no changes.

In the optimization of the entropy inequality (Section 4.1.4), no modifications are needed for Case 1. For Case 2, the differences begin in (4.40), where $F_k$ must be replaced by

$$G_k = \rho \left[ \frac{\bar{\eta}_{C_k}/l}{\alpha(\bar{\eta}_{C_k}/l)} - \frac{\rho}{\alpha(\rho)} \right] - \rho\gamma'(\rho)\left( \frac{\bar{\eta}_{C_k}}{l} - \rho \right).$$

The inequalities corresponding to (4.41) and (4.42) are obtained as for (5.9) and (5.13). After this point, the proof makes no further reference to $\tilde{c}$, and no changes are required.

5.2. *Case of small density.* In this case we have to prove the analog of Lemma 4.4 with $\rho h(\eta_x)$ in place of $c(\eta_x)$. We replace the first line of (4.75) by

$$(5.14) \qquad \nu_\Lambda^N\left[ f, \rho \sum_{x \in \Lambda} h(\eta_x) \right] = \nu_\Lambda^N\left[ f, \sum_{x \in \Lambda} \left( \rho h(\eta_x) - \frac{\rho^2}{\alpha} \right) \right].$$

The proof is identical to that of Lemma 4.4, except that the function $\varphi(\alpha, n)$ in (4.77) is now replaced by

$$\varphi(\alpha, n) = \rho \frac{n+1}{c(n+1)} - \frac{\rho^2}{\alpha}.$$

Equalities (4.77) and (4.83) hold, so no further change is needed.

## 6. Comparison between canonical and grand canonical measure: proof of Proposition 4.1.

We begin by two technical lemmas, giving uniform estimates on the canonical measure, in different regimes of density. Define $p_\Lambda^\rho(n) := \mu_\rho[\bar{\eta}_\Lambda = n]$ for $\rho > 0$, $n \in \mathbb{N}$ and $\Lambda \subset\subset \mathbb{Z}^d$. The idea is to get a Poisson approximation of $p_\Lambda^{N/|\Lambda|}(n)$ for very small values of $N/|\Lambda|$ and to use the uniform local limit theorem (see Theorem 6.1 in [8]) for the other cases.

LEMMA 6.1. *For every $N_0 \in \mathbb{N} \setminus \{0\}$, there exists a finite constant $A_0$ such that*

$$\sup_{0 < n \leq N \leq N_0} \left| p_\Lambda^{N/|\Lambda|}(n) - \frac{N^n}{n!}e^{-N} \right| \leq \frac{A_0}{|\Lambda|}$$

*for any $\Lambda \subset\subset \mathbb{Z}^d$.*

PROOF. Let $\rho := N/|\Lambda|$ and assume, without loss of generality, that $0 \in \Lambda$. Notice that

$$\mu_\rho[\bar{\eta}_\Lambda = n]$$

$$= \mu_\rho \Big[ \bar{\eta}_\Lambda = n \,\Big|\, \max_{x \in \Lambda} \eta_x \leq 1 \Big] \mu_\rho \Big[ \max_{x \in \Lambda} \eta_x \leq 1 \Big]$$

$$+ \mu_\rho \Big[ \bar{\eta}_\Lambda = n \,\Big|\, \max_{x \in \Lambda} \eta_x > 1 \Big] \mu_\rho \Big[ \max_{x \in \Lambda} \eta_x > 1 \Big].$$

We begin by proving that

$$(6.1) \qquad \mu_\rho \Big[ \max_{x \in \Lambda} \eta_x > 1 \Big] = O(|\Lambda|^{-1}),$$

uniformly in $0 < N \leq N_0$.

Indeed,

$$\mu_\rho \Big[ \max_{x \in \Lambda} \eta_x > 1 \Big] = 1 - (1 - \mu_\rho[\eta_0 > 1])^{|\Lambda|}$$

and

$$\mu_\rho[\eta_0 > 1] = \frac{1}{Z(\rho)} \sum_{k=2}^{+\infty} \frac{\alpha(\rho)^k}{c(k)!} = \frac{\alpha(\rho)^2}{Z(\rho)} \sum_{k=0}^{+\infty} \frac{\alpha(\rho)^k}{c(k+2)!}$$

$$= \frac{\alpha(\rho)^2}{Z(\rho)} \sum_{k=0}^{+\infty} \frac{c(k)!}{c(k+2)!} \frac{\alpha(\rho)^k}{c(k)!}.$$

Since $c(k)!/c(k+2)!$ is uniformly bounded, we have $\mu_\rho[\eta_0 > 1] \leq B_1 \alpha(\rho)^2 = O(|\Lambda|^{-2})$, uniformly in $0 < N \leq N_0$. Thus,

$$(1 - \mu_\rho[\eta_0 > 1])^{|\Lambda|} = (1 - O(|\Lambda|^{-2}))^{|\Lambda|} = O(|\Lambda|^{-1}),$$

which establishes (6.1).

Now let

$$\tilde{\rho} := \sum_{k=0}^{+\infty} k \mu_\rho[\eta_0 = k | \eta_0 \leq 1].$$

A trivial calculation shows that

$$\tilde{\rho} = \frac{\mu_\rho[\eta_0 \mathbf{1}(\eta_0 \leq 1)]}{\mu_\rho[\eta_0 \leq 1]} = \frac{\mu_\rho[\eta_0] - \mu_\rho[\eta_0 \mathbf{1}(\eta_0 > 1)]}{\mu_\rho[\eta_0 \leq 1]} = \frac{\rho - \mu_\rho[\eta_0 \mathbf{1}(\eta_0 > 1)]}{\mu_\rho[\eta_0 \leq 1]}.$$

Moreover,

$$\mu_\rho[\eta_0 \mathbf{1}(\eta_0 > 1)] = \frac{1}{Z(\rho)} \sum_{k=2}^{+\infty} \frac{k \alpha(\rho)^k}{c(k)!} = \frac{\alpha(\rho)^2}{Z(\rho)} \sum_{k=0}^{+\infty} \frac{(k+2) \alpha(\rho)^k}{c(k+2)!}$$

$$\leq \frac{B_2 \alpha(\rho)^2}{Z(\rho)} \sum_{k=0}^{+\infty} \frac{(k+2) \alpha(\rho)^k}{c(k)!}$$

$$= B_2 \alpha(\rho)^2 (\rho + 2) = O(|\Lambda|^{-2}),$$

and finally,

$$(6.2) \qquad \tilde{\rho} = \frac{\rho + O(|\Lambda|^{-2})}{1 + O(|\Lambda|^{-2})} = \rho + O(|\Lambda|^{-2}).$$

Observe that, for any $n \in \{0, \ldots, |\Lambda|\}$, we have

$$(6.3) \qquad \mu_\rho\Big[\bar{\eta}_\Lambda = n \,\Big|\, \max_{x \in \Lambda} \eta_x \leq 1\Big] = \binom{|\Lambda|}{n} \tilde{\rho}(1 - \tilde{\rho})^{|\Lambda|-n}.$$

This comes from the fact that the random variables $\{\eta_x : x \in \Lambda\}$, under the probability measure $\mu_\rho[\cdot \mid \max_{x \in \Lambda} \eta_x \leq 1]$, are Bernoulli independent random variables with mean $\tilde{\rho}$. The remaining part of the proof follows the classical argument of approximation of the binomial distribution with the Poisson distribution. Using (6.2) and (6.3), after some simple calculations, we get

$$\mu_\rho\Big[\bar{\eta}_\Lambda = n \,\Big|\, \max_{x \in \Lambda} \eta_x \leq 1\Big]$$

$$= \binom{|\Lambda|}{n} \tilde{\rho}(1 - \tilde{\rho})^{|\Lambda|-n}$$

$$= \frac{|\Lambda|!}{n!(|\Lambda|-n)!} \Big[\frac{N}{|\Lambda|} + O(|\Lambda|^{-2})\Big]^n \Big[1 - \frac{N}{|\Lambda|} + O(|\Lambda|^{-2})\Big]^{|\Lambda|-n}$$

$$= \frac{1}{n!}[N + O(|\Lambda|^{-1})]^n \Big[1 - \frac{N}{|\Lambda|} + O(|\Lambda|^{-2})\Big]^{|\Lambda|}$$

$$\qquad \times \frac{|\Lambda|(|\Lambda|-1)\cdot \cdots \cdot(|\Lambda|-n+1)}{|\Lambda|^n} \Big[1 - \frac{N}{|\Lambda|} + O(|\Lambda|^{-2})\Big]^{-n}$$

$$= \frac{N^n}{n!} e^{-N} + O(|\Lambda|^{-1})$$

uniformly in $0 < N \leq N_0$ and $0 \leq n < N$. This proves that there exist positive constants $v_1$ and $B_3$ such that if $\Lambda \subset\subset \mathbb{Z}^d$ is such that $|\Lambda| > v_1$, then

$$\Big| p_\Lambda^{N/|\Lambda|}(n) - \frac{N^n}{n!} e^{-N}\Big| \leq \frac{B_3}{|\Lambda|}$$

uniformly in $N \in \mathbb{N} \setminus \{0\}$ with $N \leq N_0$ and $n \in \mathbb{N}$ with $n \leq N$. The general case follows easily because the set of $n \in \mathbb{N}$, $N \in \mathbb{N} \setminus \{0\}$ and $\Lambda \subset\subset \mathbb{Z}^d$ such that $n \leq N \leq N_0$, $|\Lambda| \leq v_1$ and $0 \in \Lambda$ is finite. $\square$

PROPOSITION 6.2. *For any $\rho_0 > 0$, there exist finite positive constants $A_0$, $n_0$ and $v_0$ such* 

1.

$$(6.4) \qquad \sup_{n \in \mathbb{N}} \Big| \sqrt{\sigma^2(\rho)|\Lambda|} p_\Lambda^\rho(n) - \frac{1}{\sqrt{2\pi}} e^{-(n-\rho|\Lambda|)^2/(2\sigma^2(\rho)|\Lambda|)} \Big| \leq \frac{A_0}{\sqrt{\sigma^2(\rho)|\Lambda|}}$$

*for any $\rho \leq \rho_0$ and any $\Lambda \subset\subset \mathbb{Z}^d$ such that $\sigma^2(\rho)|\Lambda| \geq n_0$.*

2.

$$(6.5) \qquad \sup_{\substack{\rho > \rho_0 \\ n \in \mathbb{N}}} \left| \sqrt{\sigma^2(\rho)|\Lambda|} p_\Lambda^\rho(n) - \frac{1}{\sqrt{2\pi}} e^{-(n-\rho|\Lambda|)^2/(2\sigma^2(\rho)|\Lambda|)} \right| \leq \frac{A_0}{\sqrt{|\Lambda|}}$$

for any $\Lambda \subset\subset \mathbb{Z}^d$ such that $|\Lambda| \geq v_0$.

PROOF. This is a special case of the local limit theorem for $\mu_\rho$ (see Theorem 6.1 in [8]). □

PROOF OF PROPOSITION 4.1. Assume that $n := \nu_\Lambda \sigma \Lambda' \in \{0, \ldots, N\}$, and observe that

$$\nu_\Lambda^N[\eta_{\Lambda'} = \sigma_{\Lambda'}] = \frac{\mu_{N/|\Lambda|}[\eta_{\Lambda'} = \sigma_{\Lambda'}, \bar{\eta}_\Lambda = N]}{\mu_{N/|\Lambda|}[\bar{\eta}_\Lambda = N]}$$

$$= \frac{\mu_{N/|\Lambda|}[\eta_{\Lambda'} = \sigma_{\Lambda'}, \bar{\eta}_{\Lambda \setminus \Lambda'} = N - n]}{\mu_{N/|\Lambda|}[\bar{\eta}_\Lambda = N]}$$

$$= \mu_{N/|\Lambda|}[\eta_{\Lambda'} = \sigma_{\Lambda'}] \cdot \frac{\mu_{N/|\Lambda|}[\bar{\eta}_{\Lambda \setminus \Lambda'} = N - n]}{\mu_{N/|\Lambda|}[\bar{\eta}_\Lambda = N]}.$$

Thus, we have to bound from above the ratio

$$\frac{\mu_{N/|\Lambda|}[\bar{\eta}_{\Lambda \setminus \Lambda'} = N - n]}{\mu_{N/|\Lambda|}[\bar{\eta}_\Lambda = N]} = \frac{p_{\Lambda \setminus \Lambda'}^{N/|\Lambda|}(N - n)}{p_\Lambda^{N/|\Lambda|}(N)}.$$

Fix $\rho_0 > 0$ and consider three different cases.

*Small density case.* There exists $N_0 > 0$ such that

$$\sup_{\substack{N_0 \leq N \leq \rho_0|\Lambda| \\ n \in \mathbb{Z}, \Lambda \subset\subset \mathbb{Z}^d}} \frac{p_{\Lambda \setminus \Lambda'}^{N/|\Lambda|}(N - n)}{p_\Lambda^{N/|\Lambda|}(N)} < +\infty.$$

*Proof of small density case.* Assume $N/|\Lambda| \leq \rho_0$. By point 1 of Proposition 6.2, there exist positive constants $B_0$ and $n_0$ such that

$$\left| \sqrt{\sigma^2(N/|\Lambda|)|\Lambda \setminus \Lambda'|} p_{\Lambda \setminus \Lambda'}^{N/|\Lambda|}(N - n) \right.$$

$$\left. - \frac{1}{\sqrt{2\pi}} e^{-(N-n-(N/|\Lambda|)|\Lambda \setminus \Lambda'|)^2/(2\sigma^2(N/|\Lambda|)|\Lambda \setminus \Lambda'|)} \right|$$

$$\leq \frac{B_0}{\sqrt{\sigma^2(N/|\Lambda|)|\Lambda \setminus \Lambda'|}}$$

uniformly in $N/|\Lambda| \leq \rho_0$, $\sigma^2(N/|\Lambda|)|\Lambda \setminus \Lambda'| \geq n_0$ and $n \in \mathbb{Z}$. Moreover, by (2.5), there exist $B_1 > 0$ such that

$$\sigma^2(N/|\Lambda|)|\Lambda \setminus \Lambda'| \geq B_1^{-1} N \frac{|\Lambda \setminus \Lambda'|}{|\Lambda|} \geq B_1^{-1}(1-\delta_0)N \geq n_0,$$

for any $N \geq B_1 n_0/(1-\delta_0)$. Similarly, for any $N \geq 8\pi B_0^2 B_1/(1-\delta_0)$, we have

$$\sigma^2(N/|\Lambda|)|\Lambda \setminus \Lambda'| \geq B_1^{-1}(1-\delta_0)N \geq 8\pi B_0^2.$$

Thus, for any $N \geq N_0 := \lceil [n_0 \vee (8\pi B_0^2)]B_1(1-\delta_0)^{-1} \rceil$ such that $N \leq \rho_0 |\Lambda|$, we get

$$
\begin{aligned}
\sqrt{\sigma^2(N/|\Lambda|)|\Lambda \setminus \Lambda'|} \, p_{\Lambda \setminus \Lambda'}^{N/|\Lambda|}(N-n) &\leq \frac{1}{\sqrt{2\pi}} + \frac{B_0}{\sqrt{\sigma^2(N/|\Lambda|)|\Lambda \setminus \Lambda'|}} \\
&\leq \frac{1}{\sqrt{2\pi}} + \frac{1}{\sqrt{8\pi}} = \frac{3}{2\sqrt{2\pi}},
\end{aligned}
$$
(6.6)

uniformly in $n \in \mathbb{N}$. Similarly, since

$$\left| \sqrt{\sigma^2(N/|\Lambda|)|\Lambda|} \, p_\Lambda^{N/|\Lambda|}(N) - \frac{1}{\sqrt{2\pi}} \right| \leq \frac{B_0}{\sqrt{\sigma^2(N/|\Lambda|)|\Lambda|}}$$

uniformly in $N/|\Lambda| \leq \rho_0$, $\sigma^2(N/|\Lambda|)|\Lambda| \geq n_0$ and $n \in \mathbb{Z}$, by choosing $N \geq N_0$, again we get

$$
\begin{aligned}
\sqrt{\sigma^2(N/|\Lambda|)|\Lambda|} \, p_\Lambda^{N/|\Lambda|}(N) &\geq \frac{1}{\sqrt{2\pi}} - \frac{B_0}{\sqrt{\sigma^2(N/|\Lambda|)|\Lambda|}} \\
&\geq \frac{1}{\sqrt{2\pi}} - \frac{1}{\sqrt{8\pi}} = \frac{1}{2\sqrt{2\pi}}
\end{aligned}
$$
(6.7)

for any $N \geq N_0$ such that $N \leq \rho_0 |\Lambda|$ and for any $n \in \mathbb{Z}$. By (6.6) and (6.7), we get

$$\frac{p_{\Lambda \setminus \Lambda'}^{N/|\Lambda|}(N-n)}{p_\Lambda^{N/|\Lambda|}(N)} \leq 3\sqrt{\frac{|\Lambda|}{|\Lambda \setminus \Lambda'|}} \leq \frac{3}{\sqrt{1-\delta_0}}$$

for any $n \in \mathbb{Z}$, $N_0 \leq N \leq \rho_0 |\Lambda|$.

*Very small density case.* For any fixed $N_0 > 0$,

$$\sup_{\substack{0 < N \leq N_0 \\ n \in \mathbb{Z}, \Lambda \subset\subset \mathbb{Z}^d}} \frac{p_{\Lambda \setminus \Lambda'}^{N/|\Lambda|}(N-n)}{p_\Lambda^{N/|\Lambda|}(N)} < +\infty.$$

*Proof of very small density case.* In this case we use the trivial bound $p_{\Lambda\setminus\Lambda'}^{N/|\Lambda|}(N-n) \leq 1$ on the numerator and for the denominator, we use Lemma 6.1, which implies that there exists a positive constant $B_2$ such that, for any $0 < N \leq N_0$,

$$p_\Lambda^{N/|\Lambda|}(N) \geq \frac{N^N}{N!}e^{-N} - \frac{B_2}{|\Lambda|}$$

for any $\Lambda \subset\subset \mathbb{Z}^d$. It follows that

$$p_\Lambda^{N/|\Lambda|}(N) \geq \sup_{0<N\leq N_0} \frac{N^N}{N!}e^{-N} - \frac{B_2}{|\Lambda|} \geq \frac{1}{2} \sup_{0<N\leq N_0} \frac{N^N}{N!}e^{-N}$$

for any $0 < N \leq N_0$, $\Lambda \subset\subset \mathbb{Z}^d$ with

$$|\Lambda| \geq v_1 := \frac{2B_2}{\sup_{0<N\leq N_0}(N^N/N!)e^{-N}}.$$

This proves that

$$\sup_{\substack{0<N\leq N_0 \\ n\in\mathbb{Z}, |\Lambda|\geq v_1}} \frac{p_{\Lambda\setminus\Lambda'}^{N/|\Lambda|}(N-n)}{p_\Lambda^{N/|\Lambda|}(N)} < +\infty,$$

the general case follows trivially because the set of $N \in \mathbb{N} \setminus \{0\}$ and $0 \in \Lambda \subset\subset \mathbb{Z}^d$ such that $N \leq N_0$ such that $|\Lambda| < v_1$ is finite.

*Normal and large density case.* There exists $v_0 > 0$ such that

$$\sup_{\substack{N>\rho_0|\Lambda|, n\in\mathbb{Z} \\ \Lambda\subset\subset\mathbb{Z}^d, |\Lambda|>v_0}} \frac{p_{\Lambda\setminus\Lambda'}^{N/|\Lambda|}(N-n)}{p_\Lambda^{N/|\Lambda|}(N)} < +\infty.$$

*Proof of normal and large density case.* Assume that $N/|\Lambda| > \rho_0$. Then, by point 2 of Proposition 6.2, there exist positive constants $B_3$ and $v_2$, the latter depending on $\delta_0$, such that, for any $\Lambda' \subset \Lambda \subset\subset \mathbb{Z}^d$ with $|\Lambda| \geq v_2$ and $|\Lambda'|/|\Lambda| \leq \delta_0$, we have

$$\left| \sqrt{\sigma^2(N/|\Lambda|)|\Lambda \setminus \Lambda'|} p_{\Lambda\setminus\Lambda'}^{N/|\Lambda|}(N-n) \right.$$

$$\left. - \frac{1}{\sqrt{2\pi}} \exp\left[ -\frac{(N-n-(N/|\Lambda|)|\Lambda\setminus\Lambda'|)^2}{2\sigma^2(N/|\Lambda|)|\Lambda\setminus\Lambda'|} \right] \right| \leq \frac{B_3}{\sqrt{|\Lambda\setminus\Lambda'|}}$$

and

$$\left| \sqrt{\sigma^2(N/|\Lambda|)|\Lambda|} p_\Lambda^{N/|\Lambda|}(N) - \frac{1}{\sqrt{2\pi}} \right| \leq \frac{B_3}{\sqrt{|\Lambda|}}$$

uniformly in $N/|\Lambda| > \rho_0$, and $n \in \mathbb{Z}$. Now take

$$v_0 := v_2 \vee \frac{8\pi B_3^2}{1 - \delta_0}.$$

Then for any $\Lambda$ such that $|\Lambda| \geq v_0$, we have

$$\sqrt{\sigma^2(N/|\Lambda|)|\Lambda \setminus \Lambda'|} p_{\Lambda \setminus \Lambda'}^{N/|\Lambda|}(N - n)$$

$$\leq \frac{1}{\sqrt{2\pi}} + \frac{B_3}{\sqrt{|\Lambda \setminus \Lambda'|}} \leq \frac{1}{\sqrt{2\pi}} + \frac{B_3}{\sqrt{(1 - \delta_0)|\Lambda|}}$$

$$\leq \frac{1}{\sqrt{2\pi}} + \frac{1}{2\sqrt{2\pi}} = \frac{3}{2\sqrt{2\pi}}$$

and

$$\sqrt{\sigma^2(N/|\Lambda|)|\Lambda|} p_\Lambda^{N/|\Lambda|}(N) \geq \frac{1}{\sqrt{2\pi}} - \frac{B_3}{\sqrt{|\Lambda|}} \geq \frac{1}{\sqrt{2\pi}} - \frac{1}{2\sqrt{2\pi}} = \frac{1}{2\sqrt{2\pi}},$$

uniformly in $N/|\Lambda| > \rho_0$, and $n \in \mathbb{Z}$. This implies

$$\frac{p_{\Lambda \setminus \Lambda'}^{N/|\Lambda|}(N - n)}{p_\Lambda^{N/|\Lambda|}(N)} \leq 3\sqrt{\frac{|\Lambda|}{|\Lambda \setminus \Lambda'|}} \leq \frac{3}{\sqrt{1 - \delta_0}}$$

for any $n \in \mathbb{Z}$, $N > \rho_0|\Lambda|$, $\Lambda' \subset \Lambda \subset\subset \mathbb{Z}^d$ with $|\Lambda| \geq v_0$ and $|\Lambda'|/|\Lambda| \leq \delta_0$.

This completes the proof of (4.6). $\quad\square$

**Acknowledgments.** The authors wish to thank Lorenzo Bertini and Claudio Landim for very useful and interesting conversations on this work and Anna Maria Paganoni for her help in the first stage of the work.

DIPARTIMENTO DI MATEMATICA PURA E APPLICATA
UNIVERSITÀ DI PADOVA
VIA BELZONI 7
35131 PADOVA
ITALY
E-MAIL: daipra@math.unipd.it

DIPARTIMENTO DI MATEMATICA
POLITECNICO DI MILANO
PIAZZA L. DA VINCI 32
20133 MILANO
ITALY
E-MAIL: gustavo.posta@polimi.it